\documentclass{article} 
\usepackage{amsmath,amssymb,amsthm} 
\usepackage[OT2,T1]{fontenc}
\usepackage[utf8]{inputenc}
\usepackage{enumitem}
\usepackage{graphicx}
\usepackage{geometry}
\usepackage{hyperref}
\usepackage{tikz}
\usetikzlibrary{cd}
\usepackage{tikz-cd}
\usepackage{xcolor}
\usepackage{ stmaryrd }
\usepackage{mathrsfs}
\usepackage{titling}
\usepackage{fullpage}
\usepackage{mathtools}
\usepackage[english]{babel}
\usepackage{setspace}
\usepackage[normalem]{ulem}
\usepackage{mathdots}
\usepackage{microtype}

\usepackage{arydshln}

\usetikzlibrary{arrows,positioning,decorations.pathmorphing,
	decorations.markings
}
\tikzset{inner sep=0pt,
	root/.style={circle,draw,minimum size=7pt,thick},
	fatroot/.style={circle,draw,minimum size=10pt,thick},
	short root/.style={circle,fill,minimum size=7pt},
	doublearrow/.style={postaction={decorate},
		decoration={markings,mark=at position .7
			with {\arrow{angle 60}}},double distance=3pt,thick}
}

\DeclareUnicodeCharacter{00A0}{ }

\newtheorem{proposition}{Proposition}[section]
\newtheorem{definition}[proposition]{Definition}

\newtheorem{theorem}[proposition]{Theorem}
\newtheorem{lemma}[proposition]{Lemma}
\newtheorem{corollary}[proposition]{Corollary}

\numberwithin{equation}{section}

\newcommand{\G}{\mathbb{G}}
\DeclareMathOperator{\GL}{GL}

\DeclareMathOperator{\SL}{SL}

\DeclareMathOperator{\SO}{SO}

\DeclareMathOperator{\Hom}{Hom}
\DeclareMathOperator{\Mat}{Mat}

\DeclareMathOperator{\End}{End}

\DeclareMathOperator{\vol}{vol}
\DeclareMathOperator{\Vol}{vol}

\DeclareMathOperator{\tr}{tr}
\DeclareMathOperator{\rank}{rank}

\DeclareMathOperator{\Spec}{Spec}

\DeclareMathOperator{\Frac}{Frac}

\DeclareMathOperator{\Pic}{Pic}

\newcommand{\A}{\mathbb{A}}
\renewcommand{\P}{\mathbb{P}}

\newcommand{\cO}{\mathcal{O}}
\newcommand{\cR}{\mathcal{R}}

\newcommand{\R}{\mathbb{R}}
\newcommand{\Q}{\mathbb{Q}}
\newcommand{\bC}{\mathbb{C}}
\newcommand{\Z}{\mathbb{Z}}
\newcommand{\F}{\mathbb{F}}

\DeclareMathOperator{\Sel}{Sel}

\DeclareSymbolFont{cyrletters}{OT2}{wncyr}{m}{n}
\DeclareMathSymbol{\Sha}{\mathalpha}{cyrletters}{"58}

\newcommand{\extp}{\@ifnextchar^\@extp{\@extp^{\,}}}
\def\@extp^#1{\mathop{\bigwedge\nolimits^{\!#1}}}

\DeclareUnicodeCharacter{00A0}{ }
\newcommand{\define}[1]{{\fontfamily{cmss}\selectfont{#1}}}

\newcommand{\height}{\mathrm{Ht}}

\DeclareMathOperator{\covol}{covol}
\DeclareMathOperator{\disc}{disc}

\newcommand{\Siegel}{\mathfrak{S}}
\DeclareMathOperator{\irr}{irr}

\title{100\% of odd hyperelliptic Jacobians have no rational points of small height}
\author{Jef Laga and Jack A. Thorne}

\begin{document}

\maketitle

\begin{abstract}
    We study the family of hyperelliptic curves over $\Q$ given by the equation 
    \[ y^2 = x^{2g+1} + c_2 x^{2g-1} + \dots + c_{2g+1}, \]
    relating the height of rational points on their Jacobians to the reduction theory of the representation of $\SO_{2g+1}$ on self-adjoint $(2g+1)\times (2g+1)$-matrices. Using this theory, we show that in a density 1 subset, the Jacobians of these curves have no nontrivial rational points of small height.
\end{abstract}

\tableofcontents

\section{Introduction}

\paragraph{Context.} Let $g \geq 1$ be an integer. In this paper, we study the arithmetic of the family of curves $C_f$ given by the equation
\begin{equation}\label{eqn_intro_family_of_curves} C_f^0 : y^2 = f(x) 
\end{equation}
for polynomials $f(x) = x^{2g+1} + c_2 x^{2g-1} + \dots + c_{2g+1} \in \Z[x]$ of nonzero discriminant. More precisely, such an equation defines a smooth affine curve $C_f^0$ over $\Q$, and its smooth projective completion $C_f$, which is a hyperelliptic curve of genus $g$, can be obtained by adding a single rational point $P_\infty$. 

This family of curves has been the subject of a number of remarkable conjectures and theorems. Poonen and Rains \cite{PoonenRainsRandom} gave a conjectural probability distribution for the isomorphism class of the $n$-Selmer group $\Sel_n(J_f)$ of the Jacobian variety $J_f = \Pic^0_{C_f}$; this would imply in particular that the rank of $J_f(\Q)$ is either $0$ or $1$ with equal probability (and is $> 1$ with probability $0$). Here and below `probability' should be taken in the sense of natural density, with respect to the height 
\[ \height(f) = \max_i |c_i|^{1/i} \]
on the family of polynomials $f$ defining the curves $C_f$. 

Bhargava and Gross \cite{BhargavaGross} proved unconditionally that the average size of $\Sel_2(J_f)$ equals the value predicted by Poonen and Rains, namely 3. Their work was used by Poonen and Stoll \cite{PoonenStollMost} to show that when $g \geq 3$, a positive proportion of the curves in the family (\ref{eqn_intro_family_of_curves}) have no nontrivial rational points (i.e.\ satisfy $C_f(\Q) = \{ P_\infty \}$), and that this proportion rapidly approaches $1$ as $g \to \infty$. 

\paragraph{Results of this paper.} In this paper, we consider the rational points of the Jacobian variety $J_f$, and show that, in a density 1 family, nontrivial rational points must have large height, in a precise sense. 
Let $\mathcal{F}(X)$ be the set of polynomials $f(x) = x^{2g+1} + c_2 x^{2g-1} + \dots + c_{2g+1} \in \Z[x]$ of nonzero discriminant and height $<X$.
Our main theorem is as follows:
\begin{theorem}[Theorem \ref{thm_main_theorem_on_heights}]\label{thm_intro_main_theorem}
    Let $\epsilon > 0$. Then we have 
    \[ \lim_{X \to \infty} \frac{ \# \{ f \in \mathcal{F}(X)\colon  \exists P \in J_f(\Q) - \{ 0 \}, h^\dagger(P) \leq (g - \epsilon) \log \height(f) \} }{ \# \mathcal{F}(X)} = 0. \]
\end{theorem}
Here $h^\dagger(P)$ is a naive height on $J_f(\Q)$ that can be defined concretely as follows. Any divisor class $P \in J_f(\Q)$ contains a unique representative divisor $D$ which is reduced, in the sense that we can express $D = \sum_{i=1}^m P_i - m P_\infty$ for some points $P_i \in C_f^0(\overline{\Q})$ such that $P_i \neq \iota(P_j)$ if $i \neq j$ (where $\iota : C_f \to C_f$ is the hyperelliptic involution), and $m \leq g$. Then $h^\dagger(P)$ is the logarithmic Weil height of the polynomial $U(x) = \prod_i (x - x(P_i)) \in \Q[x]$. (This polynomial is part of the \emph{Mumford representation} of the divisor class $P$.)

If $g = 1$, so $C_f$ is an elliptic curve with origin $P_\infty$, this is the usual logarithmic Weil height of the $x$-coordinate of the rational point $P \in C_f(\Q)$; in general, $h^\dagger$ should be seen as a lower bound for the naive height associated to a projective embedding of the Jacobian (see Proposition \ref{prop_lower_bound_for_usual_naive_height}).

\paragraph{Methods.} We now describe our approach to the proof of Theorem \ref{thm_intro_main_theorem}. Bhargava and Gross proved their theorem on the average size of $\Sel_2(J_f)$ by relating the theory of 2-descent on the Jacobian $J_f$ to the arithmetic invariant theory of the representation of the split group $G = \SO_{2g+1}$ on the space $V$ of traceless self-adjoint linear operators on the $(2g+1)$-dimensional standard representation $W$. They showed that, if $f(x) = x^{2g+1} + c_2 x^{2g-1} + \dots + c_{2g+1} \in \Z[x]$ has nonzero discriminant, then there exists an injective map $J_f(\Q)/2J_f(\Q) \rightarrow G(\Q)\backslash V_f(\Q)$, where the target is the set of $G(\Q)$-orbits on operators $T\in V(\Q)$ of characteristic polynomial $f(x)$.
Roughly speaking, they use this to show that $\# G(\Z) \backslash V_f(\Z)$ is an upper bound for the size of $\Sel_2(J_f)$, and that a careful study of the geometry of numbers of the pair $(G, V)$ can be used to compute the exact average of $\# \Sel_2(J_f)$. 

Here we show that the reduction covariant, defined in \cite{Thorne-reduction} for an arbitrary stably graded Lie algebra, can be used to associate a lattice to a rational point $P \in J_f(\Q)$. The reduction covariant is a $G(\R)$-equivariant map 
\[ \mathcal{R} : V^{\Delta \neq 0}(\R) \to X_G \]
\[ T \mapsto H_T\]
from the set of traceless self-adjoint linear operators $T \in V(\R)$ which are regular semisimple (equivalently, of nonzero discriminant) to the symmetric space $X_G$ of the group $G(\R)$. The space $X_G$ can be interpreted as the set of inner products on $W(\R)$ which are compatible with the symmetric bilinear form defining the group $G = \SO(W)$, in a suitable sense. In particular, this construction implies that to any $T \in V(\Z)$ of nonzero discriminant, we can associate an inner product $\langle \cdot, \cdot \rangle_{H_T}$ on $W(\R)$, and therefore a lattice, i.e.\ the data of the free abelian group $W(\Z)$ together with inner product $\langle \cdot, \cdot \rangle_{H_T}$ on the real vector space $W(\R)$; and the isomorphism class of this lattice depends only on the $G(\Z)$-orbit of $T$. The main insights of this paper can now be summarised as follows:
\begin{itemize}
    \item First, the reduction covariants of irreducible orbits in $G(\Z) \backslash V(\Z)^{\Delta \neq 0}$ are equidistributed, with respect to the natural measure on $G(\Z) \backslash X_G$. (A vector of nonzero discriminant is said to be irreducible if it does not lie in the rational orbit associated to the identity $0 \in J_f(\Q)$.)
    \item Second, if $f(x) = x^{2g+1}  + c_2 x^{2g-1} + \dots + c_{2g+1} \in \Z[x]$ is a polynomial of nonzero discriminant, then we construct a canonical lifting of the map $J_f(\Q) \to G(\Q) \backslash V_f(\Q)$ to a map $J_f(\Q) \to G(\Z) \backslash V_f(\Z)$.
    Composing this map with the reduction covariant, we therefore get a map 
    \begin{equation}
        J_f(\Q) \rightarrow G(\Z) \backslash X_G.
    \end{equation}
    In other words, we assign a lattice to any rational point $P \in J_f(\Q)$. There is a subtle relation between the behaviour of the point $P$ and the associated lattice.
\end{itemize}
As an example of this relation, we show that if the Mumford representation of the point $P$ gives rise to a polynomial $U(x) \in \Z[x]$, then the associated lattice contains a vector of length
\begin{equation}\label{eqn_intro_height_norm_formula} \sum_{i=1}^{2g+1} \frac{| U(\omega_i) |}{ |f'(\omega_i) |}, 
\end{equation}
where $\omega_1, \dots, \omega_{2g+1}$ are the complex roots of $f(x)$ (see Theorem \ref{thm_good_global_lattice}). A particular consequence is that if $P \in J_f(\Q)$ has \emph{small} height (in the sense of Theorem \ref{thm_intro_main_theorem}), then the associated lattice must contain a \emph{short} vector, and therefore lie `far out in the cusp' of the locally symmetric space $G(\Z) \backslash X_G$; but this can happen only rarely, because of the equidistribution of reduction covariants. This tension is what underlies the proof of Theorem \ref{thm_intro_main_theorem}. 

\paragraph{Relation to existing results.} We now discuss the relation of Theorem \ref{thm_intro_main_theorem} to other results and conjectures about heights, focusing first on the case $g=1$. In this case, the most significant conjecture is Lang's conjecture, which asserts a uniform lower bound $\hat{h}(P) \geq c \log | \Delta_E |$ for the canonical height of nontorsion points $P$ on elliptic curves $E$ over $\Q$ of minimal discriminant $\Delta_E$ \cite[p.\ 92]{LangDiophantineAnalysis}. 

Hindry and Silverman showed, using the theory of local heights, that Lang's conjecture follows from Szpiro's conjecture, i.e.\ the uniform boundedness of the Szpiro ratio \cite{HindrySilverman}. Fouvry, Nair, and Tenenbaum showed that the Szpiro ratio can be bounded in a density 1 family \cite{FouvryNairTenenbaumSzpiro}. These results together imply that there is a constant $c > 0$ and a density 1 family of polynomials $f(x)$ for which Lang's conjecture holds (with this value of $c$).

In fact, if one wants `density 1' results, one can do better: using the theory of local heights, it is relatively easy to show (see Proposition \ref{prop_statistical_Lang} below) that for any $\epsilon > 0$, there is a density 1 family in which Lang's conjecture holds in the form $\hat{h}(P) \geq \left(\frac{1}{2}-\epsilon
\right) \log \height(f)$, without using any information on naive heights. Le Boudec \cite{LeBoudec-statistical} has combined this with a nontrivial bound for the number of points of small naive height, proved using analytic techniques, to give a density 1 bound of the form $\hat{h}(P) \geq (\frac{7}{4} - \epsilon) \log \height(f)$ for any $\epsilon > 0$. 

Unlike Lang's conjecture, a lower bound $h(P) \gg \log \height(f)$ for the naive height of $P$ can hold only in a density 1 family, and not for every elliptic curve, as there are infinitely many examples of pairs $(E, P)$ where $P \in E(\Q)$ is a rational point of naive height 0.  For this reason we view a result such as Theorem \ref{thm_intro_main_theorem} as related, but far from equivalent, to possible generalisations of Lang's conjecture to abelian varieties of dimension $g > 1$, such as the one given by Silverman \cite{SilvermanHeightsAbelianVarieties}. 

In the context of higher genus curves, it is frequently useful to divide possible rational points into `small' and `large' height subsets (see e.g.\ \cite{gao2021recent}). Theorem \ref{thm_intro_main_theorem} can be seen as asserting that the `small height' subset of $J_f(\Q)$ is almost always trivial (i.e. reduced to $\{ 0 \}$). Our method gives even stronger lower bounds for the heights of rational points on the subsets of $J_f(\Q)$ whose Mumford degree (i.e.\ degree of the associated polynomial $U$) is strictly less than $g$. For example, we can show:
\begin{theorem}\label{thm_intro_degree1}
    Let $\epsilon > 0$. Then we have 
    \[ \lim_{X \to \infty} \frac{ \# \{ f \in \mathcal{F}(X)\colon  \exists P = (\alpha, \beta) \in C^0_f(\Q), h(\alpha) \leq (2g - 1 - \epsilon)\log \height(f) \} }{ \# \mathcal{F}(X)} = 0. \]
\end{theorem}
There is an interesting tension between this result and the conjecture that there are constants $c_g, d_g > 0$ such that for any $f \in \mathcal{F}(X)$ and $P = (\alpha, \beta) \in C_f^0(\Q)$, we have
\[ h(\alpha) \leq c_g \log \height(f) + d_g \]
(as shown in \cite{IhHeightUniformity}, this would follow from  Vojta's conjecture), especially in light of the expectation (see \cite[Remark 10.11]{PoonenStollMost}) that $100 \%$ of curves in the family (\ref{eqn_intro_family_of_curves}) have no rational points other than $P_\infty$.

\paragraph{Structure of this paper.} In \S \ref{sec_preliminaries} we recall the precise definition and basic structures associated to the representation $(G, V)$ of the rank $2g+1$ special orthogonal group. We also recall from \cite{Thorne-reduction} the definition and properties of the reduction covariant $\mathcal{R}$ in this case and explain how it may be computed. 

In \S \ref{section: equidistribution reduction covariant} we prove the equidistribution of the reduction covariants of irreducible orbits in $G(\Z) \backslash V(\Z)^{\Delta \neq 0}$. The key point here is to show that the definition of $\mathcal{R}$ meshes well with the approach to the geometry-of-numbers used in \cite{BhargavaGross}.

In \S \ref{sec_classificationconstructionorbits} we establish the critical results that allow us to relate the lattices determined by vectors $T \in V(\Z)^{\Delta \neq 0}$ to heights. In particular, we give a precise refinement (Theorem \ref{thm_good_global_lattice}) of \cite[Proposition 19]{BhargavaGross}, which constructs integral orbit representatives for rational orbits associated to rational points of the Jacobian $J_f$. We also make a connection with the reduction theory of Stuhler and Grayson \cite{Stuhler-reduction, Grayson-reduction}, which associates to any lattice its Harder--Narasimhan polygon (i.e.\ what is called the canonical plot in \cite{Grayson-reduction}, encoding the data of the minimal covolumes of subgroups of each possible rank of $W(\Z)$). We use this polygon when analyzing the heights of points in $2 J_f(\Q) - \{ 0 \}$. 

In the short \S \ref{sec_density_1_family}, we construct a density 1 family of polynomials $f(x) = x^{2g+1} + c_2 x^{2g-1} + \dots + c_{2g+1} \in \Z[x]$ which has good (or `good enough') properties. Polynomials whose discriminants are divisible by large prime powers give rise to curves with poorly behaved height functions; we use the methods of \cite{Bhargava-squarefree} to show that these do not occur too often. Similarly, polynomials with complex roots that  are too close together make it hard to relate quantities such as (\ref{eqn_intro_height_norm_formula}) to the heights of interest. A more elementary argument suffices to show that these also occur rarely.

Finally, in \S \ref{sec_lower_bounds_for_height_functions} we bring all of these ingredients together, following the sketch outlined above, in order to prove Theorem \ref{thm_intro_main_theorem}.

\section{Preliminaries}\label{sec_preliminaries}

\subsection{Quadratic spaces}\label{subsec_quadspaces}

Let $k$ be a field, and let $N \geq 1$ be an integer. In this paper, we call a \define{quadratic space} over $k$ of rank $N$ a $k$-vector space $W_0$ of dimension $N$, equipped with a nondegenerate symmetric bilinear form $(\cdot, \cdot)_{W_0} : W_0 \times W_0 \to k$. Its \define{discriminant} is defined to be the element $(-1)^{N(N-1)/2} \det A \text{ mod }(k^\times)^2 \in k^\times / (k^\times)^2$, where $A$ is the Gram matrix of $(\cdot, \cdot)_{W_0}$ with respect to any $k$-basis of $W_0$. A linear endomorphism $T_0 : W_0 \to W_0$ is said to be \define{self-adjoint} if it satisfies the formula $(T_0 \cdot, \cdot)_{W_0} = (\cdot, T_0 \cdot)_{W_0}$. 

We will generally work with a fixed quadratic space $W$ of rank $N = 2g +1$ and its associated space $V$ of traceless, self-adjoint linear endomorphisms; these are introduced in \S \ref{subsec_basic_definitions}. However, later in the paper, a useful method to construct elements of $W$ and $V$ will be to first construct data in auxiliary quadratic spaces and show eventually that they are isomorphic to $W$.

\subsection{The representation $(G,V)$}\label{subsec_basic_definitions}

Fix an integer $g\geq 1$. 
Let $W$ be the free $\Z$-module with basis $e_{-1}, \dots,e_{-g}, e_0 , e_g,\dots,e_1$, equipped with the symmetric bilinear form $( \cdot, \cdot)_W$ satisfying $( e_i , e_j)_W = \delta_{i, -j}$ if $-g \leq i, j \leq g$.
The Gram matrix $J$ of $(\cdot, \cdot)_W$ is thus 
\begin{align*}
    \begin{pmatrix}
         & & 1 \\
         & \iddots & \\
         1 & & \\
    \end{pmatrix}.
\end{align*}
Given a linear map $A \in \End(W) \simeq \Mat_{2g+1}(\Z)$, its adjoint $A^* = J ({}^tA) J$ is the unique linear map satisfying $( Av,w )_W =  (v, A^* w)_W$ for all $v,w\in W$. (Throughout the paper we write ${}^t A$ for the transpose of a matrix $A$.) Then $A^\ast$ is obtained from $A$ by flipping along the anti-diagonal. 

Following \cite{BhargavaGross}, define the group scheme $G$ over $\Z$ by 
\begin{align*}
    G \coloneqq \SO(W) = \{g\in \SL(W) \colon g g^* = 1 \}
\end{align*}
and the representation $V$ over $\Z$ by
\begin{align*}
    V \coloneqq \{ T\in \End(W) \colon T^* = T\text{ and } \tr(T) =0 \}.
\end{align*}
The free $\Z$-module $V(\Z)$ has rank $2g^2+3g$ and $G$ acts on $V$ via $g\cdot T \coloneqq g T g^{-1} = gT g^*$.

\subsection{Invariants of $(G,V)$}\label{subsec_invariantsofGV}

If $R$ is a ring and $T\in V(R)$, the characteristic polynomial of $T$ is a polynomial of the form 
\begin{align}\label{equation_charpolyofT}
    f(x) = \det(xI- T) = x^{2g+1} + c_2 x^{2g-1} + \dots + c_{2g+1} \in R[x].
\end{align}
The coefficients of the characteristic polynomial are invariants of the $G(R)$-action, and they generate the full ring of polynomial invariants of $V$ when $R  = \Q$.
Let $B = \Spec(\Z[c_2,\dots,c_{2g+1}])$ be the affine space parametrizing monic polynomials of the form \eqref{equation_charpolyofT}.
Taking the characteristic polynomial then defines a morphism $V\rightarrow B$. This morphism is equivariant for the action of $G$ (where $G$ acts trivially on $B$) and also for the action of $\G_m$ (where $\G_m$ acts by scalar multiplication on $V$, and by $t \cdot c_i = t^i c_i$ on $B$).
Given a polynomial $f\in R[x]$ of the above form, let $V_f(R)$ be the subset of $T\in V(R)$ with characteristic polynomial $f$.

If $T\in V(R)$ has characteristic polynomial $f$, write $\Delta(T) = \Delta(f)$ for the discriminant of $f$, an element in $R$.
If $c_i$ has degree $i$, then the discriminant is a homogeneous polynomial of degree $2g(2g+1)$ in $\Z[c_2,\cdots,c_{2g+1}]$.
Let $V^s\subset V$ and $B^s\subset B$ denote the open subschemes where the discriminant is nonzero.

Given an element $T\in V(\R)$ with characteristic polynomial $f = x^{2g+1} + c_2x^{2g-1} + \cdots + c_{2g+1}$, we define the \define{height} of $T$ and $f$ to be
\begin{align*}
     \height(T) \coloneqq \height(f)\coloneqq \max_i | c_i|^{1/i}.
\end{align*}
If $f \in B(\R)$ and $t \in \R^\times$ then we have $\height(t \cdot f) = |t| \cdot \height(f)$. 
\begin{lemma}\label{lem_height_comparison} Let $f(x) = x^{2g+1} + c_2 x^{2g-1} + \dots + c_{2g+1} \in \R[x]$ be a polynomial (possibly of discriminant $0$). Define $\height_1(f) = \max_i |\omega_i|$, where $\omega_1, \dots, \omega_{2g+1}$ are the roots of $f$ in $\bC$, with multiplicity. Then we have
        \[ \frac{1}{2} \height_1(f) \leq \height(f) \leq 2^{2g+1} \height_1(f). \]
\end{lemma}
\begin{proof}
    The inequality $|c_i| \leq 2^{2g+1} \max_j |\omega_j|^i$ follows from expanding out $f(x) = \prod_i (x - \omega_i)$ and using the triangle inequality. This implies that $\height(f) \leq 2^{2g+1} \height_1(f)$. For the other inequality, we need to show that $|\omega_j| \leq 2 \height(f)$ for every $j$. We can assume that $\omega_j \neq 0$. We have
    \[ \omega_j^{2g+1} = - (c_2 \omega_j^{2g-1} + \dots + c_{2g+1}), \]
    hence 
    \[ 1 = | c_2 / \omega_j^2 + c_3 / \omega_j^3 + \dots + c_{2g+1} / \omega_j^{2g+1} | \leq \sum_{i=2}^{2g+1} | \height(f) / \omega_j |^i. \]
    Suppose for contradiction that $| \omega_j | > 2 \height(f)$, hence $| \height(f) / \omega_j | < \frac{1}{2} $. Then we can complete the geometric series and sum to obtain
    \[ 1 \leq \sum_{i=2}^\infty 2^{-i}, \]
    a contradiction. 
\end{proof}

\subsection{The distinguished orbit}\label{subsec_distinguishedorbit}

If $k/\Q$ is a field and $f = x^{2g+1} + c_2x^{2g-1} + \cdots +c_{2g+1} \in k[x]$ a separable polynomial, the set 
\begin{align*}
   \{T\in V_f(k) \colon  \text{there exists a maximal isotropic }L \leq W(k) \text{ with } T L \leq L^{\perp}\} 
\end{align*}
consists of a single $G(k)$-orbit \cite[Proposition 6]{BhargavaGross}.
We call such an orbit a \define{distinguished orbit} over $k$ (or the distinguished orbit of characteristic polynomial $f$ over $k$). If $k$ is algebraically closed, every element of $V_f(k)$ lies in the distinguished orbit, but for a general field $k$ this is not the case, see \S\ref{subsec_rationalorbits}.

\subsection{The reduction covariant}

Let $W_0$ be a quadratic space over $\R$ of rank $2g+1$.
An \define{inner product} on $W_0$ is a positive definite symmetric bilinear form $\langle\cdot, \cdot\rangle_{H_0}$ on $W_0$.
The bilinear forms $(\cdot,\cdot)_{W_0}$ and $\langle\cdot,\cdot\rangle_{H_0}$ induce isomorphisms $\varphi_{W_0}, \varphi_{H_0}\colon W_0 \rightarrow W_0^*$ via the formulae $v\mapsto (v,\cdot)_{W_0}$ and $v\mapsto \langle v,\cdot\rangle_{H_0}$ respectively. 
We say $\langle \cdot, \cdot \rangle_{H_0}$ is \define{compatible with $(\cdot,\cdot)_{W_0}$} if $\varphi_{H_0}^{-1}\varphi_{W_0}$ is an isometry with respect to the form $\langle \cdot,\cdot \rangle_{H_0}$.

\begin{proposition}\label{prop_existence_of_reduction_covariant}
    Let $T_0 : W_0 \to W_0$ be a self-adjoint linear operator which is regular semisimple (equivalently, with separable characteristic polynomial).
    Then there exists a unique inner product $\langle \cdot, \cdot \rangle_{H_{T_0}}$ on $W_0$ satisfying the following conditions:
    \begin{enumerate}
        \item $\langle \cdot, \cdot \rangle_{H_{T_0}}$ is compatible with the quadratic form $(\cdot, \cdot )_{W_0}$. 
        \item $T_0$ commutes with its $\langle \cdot, \cdot \rangle_{H_{T_0}}$-adjoint.
    \end{enumerate}
\end{proposition}
\begin{proof}
    See \cite[Theorem 2.7]{Thorne-reduction} (which proves a much more general result) and \cite[Proposition 3.1]{Thorne-reduction} (which describes what happens more explicitly for a split quadratic space, a hypothesis which is inessential here). 
    Alternatively, this follows from the description of $\langle \cdot , \cdot \rangle_{H_{T_0}}$ given below.
\end{proof}
We can describe the inner product $\langle \cdot, \cdot \rangle_{H_{T_0}}$ explicitly as follows. Let $\omega_1, \dots, \omega_{2g+1} \in \bC$ be the roots of $f(x)$ (distinct, by assumption). Then $W_0 \otimes_\R \bC$ admits a basis $b_1, \dots, b_{2g+1}$ of $T_0$-eigenvectors, where $b_i$ has eigenvalue $\omega_i$. The inner product $\langle \cdot, \cdot \rangle_{H_{T_0}}$ extends uniquely to a Hermitian inner product on $W_0 \otimes_\R \bC$. 
Since $T_0$ commutes with its $\langle \cdot , \cdot\rangle_{H_{T_0}}$-adjoint, the $b_i$ are $\langle \cdot, \cdot \rangle_{H_{T_0}}$-orthogonal and since the bilinear forms are compatible they satisfy $\langle b_i, b_i \rangle_{H_{T_0}} = | (b_i, b_i)_{W_0} |$. This description implies the following useful property:
\begin{lemma}\label{lem_constant_on_Cartan_subspace}
    Let $T_0, T_1 : W_0 \to W_0$ be self-adjoint linear operators which are regular semisimple, and which commute. Then $\langle \cdot, \cdot \rangle_{H_{T_0}} =  \langle \cdot, \cdot \rangle_{H_{T_1}}.$
\end{lemma}
\begin{proof}
    There exists a basis of $W_0$ consisting of simultaneous eigenvectors for $T_0$ and $T_1$. 
\end{proof}
Recall that $V^s \subset V$ denotes the open subscheme where $\Delta$ is nonzero. The symmetric space $X_{G}$ is, by definition, the set of inner products $\langle \cdot, \cdot \rangle_H$ on $W(\R)$ which are compatible with $J$. The group $G(\R)$ acts on $X_G$ on the left by the formula
\[ \langle v, w \rangle_{g \cdot H} = \langle g^{-1} v, g^{-1} w \rangle_{H}. \]
The assignment $T \mapsto \langle \cdot, \cdot \rangle_{H_T}$ defines a map $\cR : V^s(\R) \to X_{G}$, which we call the \define{reduction covariant}. This map is equivariant for the action of $G(\R)$, and satisfies $\cR(\lambda v) = \cR(v)$ for all $\lambda \in \R^\times$, $v \in V^s(\R)$. By passage to quotient, we obtain a map
\begin{align}\label{equation_reduction_covariant}
    \cR : G(\Z) \backslash V^s(\R) \to G(\Z) \backslash X_G.
\end{align}
The quotient $G(\Z) \backslash X_G$ is an example of a locally symmetric space and admits a canonical (up to scaling) measure $\mu_G$, induced by the $G(\R)$-invariant measure on the symmetric space $X_G$. We will describe this measure more explicitly in \S \ref{subsec_preliminariescounting}. 

\section{Equidistribution of the reduction covariant}\label{section: equidistribution reduction covariant}

The purpose of this section is to show that the reduction covariant \eqref{equation_reduction_covariant} becomes equidistributed over irreducible integral orbits of $V(\Z)$. 
After discussing some preliminaries to make this statement precise, we state and prove this equidistribution using a modification of the orbit-counting techniques of Bhargava--Gross \cite[\S 10]{BhargavaGross}.
For the purposes of our main Theorem \ref{thm_intro_main_theorem}, we will only use the consequence (Corollary \ref{corollary_polyswithsmallnormvectorrare}) that polynomials $f\in \mathcal{F}(X)$ admitting an irreducible orbit $T\in V_f(\Z)$ and vector $w\in W(\Z) - \{ 0 \}$ of small $H_T$-norm are `sparse' in $\mathcal{F}(X)$.

\subsection{Preliminaries}\label{subsec_preliminariescounting}

\paragraph{Coordinates on $G$.}

Use the basis $e_{-1},\dots,e_1$ of $W$ from \S\ref{subsec_basic_definitions} to identify $G$ with a subgroup of $\GL_{2g+1}$.
Let $H_0$ be the standard inner product on $W(\R)$ for which this basis is orthonormal.
Then $K\coloneq G(\R)\cap \SO_{H_0}(\R)$ is a maximal compact subgroup of $G(\R)$.
Let $T \leq G$ denote the subgroup of diagonal matrices contained in $G$.
We use the isomorphism $\G_m^g \rightarrow T, (t_1,\dots,t_g) \mapsto \text{diag}(t_1^{-1},\dots,t_g^{-1},1,t_g,\dots,t_1)$ to define coordinates $(t_1,\dots,t_g)$ on $T$.
Let $N \leq G$ denote the subgroup of unipotent lower triangular matrices, and let $P = T N \leq G$ denote the subgroup of lower triangular matrices. 
By the Iwasawa decomposition, the product map $N(\R)\times T(\R)^{\circ} \times K \rightarrow G(\R)$ is a diffeomorphism.

\paragraph{Fundamental set for $G(\Z)\backslash G(\R)$.}

In the remainder of \S\ref{section: equidistribution reduction covariant} we will use the language of semialgebraic sets, see \cite[Chapter 2]{BCR-realalgebraicgeometry}. We will use without further mention that semialgebraic sets are closed under finite unions, intersections, and (pre-)images under semialgebraic maps.

In this paper we call a \define{Siegel set} a subset of $G(\R)$ of the form $\omega \cdot T_c \cdot K$, where $\omega \subset N(\R)$ is a relatively compact subset, $c \in \R_{>0}$ and $T_c \coloneqq  \{(t_1,\dots,t_g) \in T(\R) \colon t_1/t_2> c, \dots,t_{g-1}/t_g > c, t_g> c\}$.
For every such Siegel set, the set of $\gamma\in G(\Z)$ with $\gamma\cdot  \Siegel \cap \Siegel \neq \varnothing$ is finite \cite[Corollaire 15.3]{Borel-introductiongroupesarithmetiques} (`Siegel property'). 
By Proposition \ref{prop_good_neighbourhoods_of_cusp} and \cite[Proposition 15.6]{Borel-introductiongroupesarithmetiques}, there exists a Siegel set $\Siegel$ with the property that $G(\Z)\cdot \Siegel = G(\R)$ \cite[\S9.2]{BhargavaGross}.
After enlarging $\Siegel$, we may assume that $\omega$, and consequently $\Siegel$, is open and semialgebraic.

The set $\Siegel$ will serve as our fundamental set for the left action of $G(\Z)$ on $G(\R)$.
A $G(\Z)$-orbit might be represented more than once in $\Siegel$, but this does not cause any problems as long as we incorporate the multiplicity function $m \colon \Siegel \rightarrow \Z_{\geq 1}$, defined by $m(x)=\#(G(\Z) \cdot x \cap \Siegel)$.
This function is bounded and has semialgebraic fibres.

\begin{lemma}\label{lemma_strongsiegelproperty}
    For every two compact subsets $A, B\subset G(\R)$, the set $\{\gamma\in G(\Z) \colon \gamma A \cap (\Siegel \cdot B) \neq \varnothing\}$ is finite.
\end{lemma}
\begin{proof}
    Since $\gamma A\cap (\Siegel \cdot B) \neq \varnothing$ if and only if $(A\cdot B^{-1})\cap \gamma^{-1} \Siegel \neq \varnothing$, it suffices to prove for every compact subset $Z\subset G(\R)$ that the set $\{\gamma \in G(\Z) \colon  Z \cap \gamma\Siegel \neq \varnothing\}$ is finite.
    Since $G(\R) = G(\Z) \cdot \Siegel$, the open subsets $\{\gamma\Siegel\}_{\gamma \in G(\Z)}$ cover $Z$.
    By compactness, there are $\gamma_1,\dots,\gamma_n$ such that $Z$ is contained in $\gamma_1\Siegel \cup \cdots \cup \gamma_n \Siegel$.
    If $\gamma\in G(\Z)$ is such that $Z \cap \gamma \Siegel \neq \varnothing$, then $\gamma \Siegel \cap \gamma_i \Siegel\neq \varnothing$ for some $i$.
    By the Siegel property, there are only finitely many $\gamma$ satisfying the latter condition for some $i$.
\end{proof}

\paragraph{Measures on $G$ and $X_G$.}
Consider the $1$-dimensional vector space of left invariant top differential forms of $G$ over $\Q$. Those forms that extend to a differential form over $\Z$ form a rank-$1$ submodule of this vector space. 
Choose a generator (which is unique up to sign) of this rank $1$ submodule.
It induces a bi-invariant Haar measure $dg$ on $G(\R)$.
We equip the maximal compact subgroup $K$ with its probability Haar measure.
We now explain how these measures also induce measures on $X_G$ and $G(\Z)\backslash X_G$.

The standard inner product $H_0$ defines an element of $X_G$ and induces a $G(\R)$-equivariant bijection $G(\R)/K \simeq X_G$.
We will use this identification without further mention.
Since the measure $dg$ on $G(\R)$ is bi-invariant, it induces measures on $G(\R)/K=X_G$ and $G(\Z)\backslash G(\R) /K=G(\Z)\backslash X_G $. 
We denote the latter measure by $\mu_G$.
It is a standard fact that $\mu_G(G(\Z)\backslash X_G)$ is finite, which follows from the fact that $\Vol(\Siegel)$ is finite. 

\paragraph{Good subsets of $X_G$.}
Consider the surjective quotient map 
$\varphi\colon \Siegel \rightarrow G(\Z)\backslash X_G$.
We call a subset $U \subset G(\Z)\backslash X_G$ \define{good} if it is relatively compact and $\varphi^{-1}(U)\subset \Siegel$ is a semialgebraic subset of $G(\R)$.
For example, the image under $\varphi$ of a relatively compact semialgebraic subset of $\Siegel$ is good.

\begin{lemma}\label{lemma_countablebasis_goodsubsets}
    There exists a countable basis of good open subsets of $G(\Z)\backslash X_G$.
\end{lemma}
\begin{proof}
    It suffices to prove that $\Siegel = \omega T_c K$ has a countable basis consisting of relatively compact and semialgebraic open subsets. Since $\Siegel$ is open in $G(\R)$, it even suffices to prove that $G(\R)$ has such a countable basis. Since $G$ is a closed subscheme of $\A^N$ for some $N\geq 1$, we may take the basis $\{ B\cap G(\R)\}$ of $G(\R)$, where $B$ runs over all open balls in $\R^N$ with centre in $\Q^N$ and rational radius.
    \end{proof}

The next two lemmas ensure that good subsets have good properties for the purposes of the geometry-of-numbers arguments.

\begin{lemma}\label{lemma_integral_siegelset}
    Let $U\subset G(\Z)\backslash X_G$ be good and let $\bar{U}\subset G(\R)$ be its preimage under the quotient map $G(\R)\rightarrow G(\Z)\backslash X_G$. 
    Then for all $h\in G(\R)$ we have
    \begin{align*}
    \int_{g \in \Siegel \cap \bar{U}h} m(g)^{-1}\,dg = \mu_G(U).
    \end{align*}
\end{lemma}
\begin{proof}
    Follows from pushing forward measures and the bi-invariance of $dg$.
\end{proof}

\begin{lemma}\label{lemma_goodsubsetimpliesregionsemialg}
    Let $U\subset G(\Z)\backslash X_G$ be good, and let $A\subset G(\R)$ be a compact semialgebraic subset. 
    Let $\bar{U}$ be the preimage of $U$ under the quotient map $G(\R) \rightarrow G(\Z)\backslash X_G$.
    Then the set $\{ (g,h) \in \Siegel \times A \colon gh \in \bar{U} \}$ is a semialgebraic subset of $G(\R)\times G(\R)$.
\end{lemma}

\begin{proof}
    Let $W = \bar{U} \cap \Siegel$.
    Since $U$ is good and $\varphi$ is proper, $W$ is semialgebraic and relatively compact.
    Since $G(\R)=G(\Z) \cdot \Siegel$, $\bar{U} = G(\Z) \cdot W$.
    By Lemma \ref{lemma_strongsiegelproperty}, only finitely many elements $\gamma_1,\dots,\gamma_n \in G(\Z)$ satisfy $\gamma W \cap (\Siegel \cdot A) \neq \varnothing$. 
    Therefore $\bar{U} \cap (\Siegel \cdot A) = ( \gamma_1 W \cup \cdots \cup \gamma_n W) \cap (\Siegel \cdot A)$.
    Since $\gamma_iW$ and $\Siegel\cdot A$ are semialgebraic, $\bar{U}\cap (\Siegel\cdot A)$ is semialgebraic.
    Since $\{(g,h) \in \Siegel \times A\colon gh\in \bar{U} \}$ is the pre-image of $\bar{U}\cap (\Siegel \cdot A)$ under the multplication map $\Siegel \times A \rightarrow G(\R)$, it is semialgebraic too.
\end{proof}

\paragraph{Fundamental sets for $G(\R)\backslash V(\R)$.}

Recall from \S\ref{subsec_invariantsofGV} that $B^s(\R)$ denotes the set of polynomials $f = x^{2g+1} +c_2x^{2g-1} + \cdots +c_{2g+1} \in \R[x]$ with nonzero discriminant. 
Let $I(m) \subset B^s(\R)$ be the subset of polynomials having exactly $2m+1$ real roots. 
The sets $I(m)$ with $0\leq m\leq g$ are the connected components of $B^s(\R)$.
Bhargava--Gross \cite[\S6.3]{BhargavaGross} have further partitioned $V^s(\R)$ into components $V^{(m,\tau)}$ indexed by integers $0\leq m\leq g$ and $1\leq \tau \leq \binom{2m+1}{m}$, whose definition we will not recall here. 
We only mention that the component $V^{(m,\tau)}$ maps onto $I(m)$ under the characteristic polynomial map $V^s(\R)\rightarrow B^s(\R)$ and $V^{(m,1)}$ consists of the distinguished orbits lying over $I(m)$.

In \cite[\S9.1]{BhargavaGross}, Bhargava--Gross construct fundamental sets $L^{(m,\tau)}$ for the action of $G(\R)$ on height $1$ elements in $V^{(m,\tau)}$.
These are bounded and semialgebraic subsets of $V^{(m,\tau)}$ with the property that the characteristic polynomial map $L^{(m,\tau)}\rightarrow \{ f\in I(m)\colon \height(f) =1 \}$ is a semialgebraic isomorphism.

\begin{lemma}
    For every pair $(m,\tau)$ as above and $x,y\in L^{(m,\tau)}$, we have $\mathcal{R}(x) = \mathcal{R}(y)$ in $X_G$.
\end{lemma}
\begin{proof}
    The explicit description of $L^{(m,\tau)}$ given in \cite[\S9.1]{BhargavaGross} shows that every two elements in $L^{(m,\tau)}$, seen as linear operators $W(\R)\rightarrow W(\R)$, commute with each other.
    We conclude using Lemma \ref{lem_constant_on_Cartan_subspace}.
\end{proof}

Since we are free to replace $L^{(m,\tau)}$ by $g\cdot L^{(m,\tau)}$ for some $g\in G(\R)$ and preserve all of its required properties, we may and do assume that the $L^{(m,\tau)}$ have been chosen so that the reduction covariant of every element in $L^{(m,\tau)}$ equals the standard inner product $H_0 \in X_G$.

\paragraph{Counting lattice points.}

We recall the following proposition \cite[Theorem 1.3]{BarroeroWidmer-lattice}, which strengthens a well-known result of Davenport \cite{Davenport-onaresultofLipschitz}.
\begin{proposition}\label{prop_countlatticepointsbarroero}
	Let $m,n\geq 1$ be integers, and let $Z\subset \R^{m+n}$ be a semialgebraic subset. 
	For $T\in \R^m$, let $Z_T = \{x\in \R^n\colon (T,x) \in Z\}$, and suppose that all such subsets $Z_T$ are bounded.
	Then
	\begin{align*}
		\#(Z_T \cap \Z^n) = \Vol(Z_T)+O(\max\{1,\Vol(Z_{T,j})\}),
	\end{align*}
	where $Z_{T,j}$ runs over all orthogonal projections of $Z_T$ to all $j$-dimensional coordinate hyperplanes $(1\leq j \leq n-1)$. 
	Moreover, the implied constant depends only on $Z$. 
 \end{proposition}

\subsection{The equidistribution theorem}

We say an element in $V(\Z)$ is \define{reducible} if it has zero discriminant or if it lies in a distinguished orbit over $\Q$, and \define{irreducible} otherwise.
For any subset $S\subset V(\Z)$, write $S^{\irr}$ for the subset of irreducible elements.
Given $X\in \R_{\geq 0}$ and a $G(\Z)$-invariant subset $S\subset V(\Z)$, let $N(S;X)$ be the number of irreducible $G(\Z)$-orbits in $S$ of height $< X$.
Let $\mathcal{F}(X)$ be the set of polynomials $f = x^{2g+1} + c_2x^{2g-1} + \cdots + c_{2g+1} \in \Z[x]$ of nonzero discriminant and of height $<X$. 
\begin{lemma}\label{lemma_F(X)_expectedsize}
    $\# \mathcal{F}(X) = 2^{2g} X^{g(2g+3)} + O(X^{g(2g+3) - 1})$.
\end{lemma}
\begin{proof}
    We need to show that there are $O(X^{g(2g+3) - 1})$ polynomials with $\height(f) < X$ and $\Delta(f) = 0$. This follows immediately from \cite[Lemma 3.1]{bhargava2014geometric}.
\end{proof}

Bhargava--Gross showed in \cite[\S10]{BhargavaGross} that the average number of irreducible integral orbits with given characteristic polynomial is finite: there exists a constant $c>0$ such that 
\begin{align}\label{equation_averageirredorbitsfinite}
N(V(\Z);X) = c\cdot \#\mathcal{F}(X)+ o(X^{g(2g+3)})
\end{align}
as $X\rightarrow +\infty$.
The main result of this section shows that the reduction covariant map $\mathcal{R}\colon G(\Z)\backslash V(\Z) \rightarrow G(\Z)\backslash X_G$ is equidistributed over irreducible integral orbits in some sense.
Recall from \S\ref{subsec_preliminariescounting} that we have equipped $G(\Z)\backslash X_G$ with a measure $\mu_G$.

\begin{theorem}\label{theorem: equidistribution reduction covariant}
    Let $U\subset G(\Z)\backslash X_G$ be a Borel-measurable subset whose boundary has measure zero.
    Then 
    \begin{align*}
        N(V(\Z)\cap \mathcal{R}^{-1}(U); X) = \frac{\mu_G(U)}{\mu_G(G(\Z)\backslash X_G) }\cdot N(V(\Z);X)+ o(X^{g(2g+3)}).
    \end{align*}
\end{theorem}

The proof of Theorem \ref{theorem: equidistribution reduction covariant} follows from a modification of the proof of Bhargava--Gross and is given below.
We first prove a local version for good subsets.
Fix a pair $(m,\tau)$ corresponding to an irreducible component $V^{(m,\tau)}\subset V^s(\R)$ as in \S\ref{subsec_preliminariescounting} and write $V(\Z)^{(m,\tau)} = V(\Z) \cap V^{(m,\tau)}$.
Let $\mathcal{J}$ be the rational nonzero constant of \cite[Proposition 34]{BhargavaGross}.
Let $c_{m,\tau} = \mu_G(G(\Z)\backslash X_G)\frac{ |\mathcal{J}| \Vol(I(m)_{<1})}{2^{m+n}}$.
Bhargava--Gross \cite[Theorem 25]{BhargavaGross} showed that 
\begin{align}\label{equation_averageorbitslocalversion}
    N(V(\Z)^{(m,\tau)};X) = c_{m,\tau} X^{g(2g+3)} + o(X^{g(2g+3)}).
\end{align}
\begin{proposition}\label{prop_equidistribution_realcomponent}
    Let $U\subset G(\Z)\backslash X_G$ be a good subset and let $S_U \coloneqq V(\Z)^{(m,\tau)} \cap \mathcal{R}^{-1}(U)$.
    Then 
    \begin{align}\label{equation_equidistribution_realcomponent}
    N(S_U;X) = \frac{\mu_G(U)}{\mu_G(G(\Z)\backslash X_G)} c_{m,\tau} X^{g(2g+3)} + o(X^{g(2g+3)}).
\end{align}
\end{proposition}
\begin{proof}
    We will use the notations and choices made in \S\ref{subsec_preliminariescounting}.
    Additionally, given a subset $A$ of $V(\R)$ or $B(\R)$, let $A_{<X}\subset A$ be the subset of elements of height $<X$.
    Write $\Lambda = \R_{>0}$, equipped with its multiplicative Haar measure $d^{\times}\lambda$, and let $\Lambda$ act on $V(\R)$ by scalar multiplication.
    Write $L  = L^{(m,\tau)}$.
    Fix a compact, semialgebraic set $G_0\subset G(\R) \times \Lambda$ of volume $1$, with nonempty interior, that satisfies $K\cdot G_0 = G_0$ and whose projection onto $\Lambda$ is contained in $[1,K_0]$ for some $K_0>1$.
    
    Given a $G(\Z)$-invariant subset $S\subset V(\Z)$, let $N^*(S;X)$ be the number of (not necessarily irreducible) $G(\Z)$-orbits in $S$ of height $<X$.
    For any subset $S\subset V(\Z)^{(m,\tau)}$, averaging over $G_0$ \cite[Equation (14)]{BhargavaGross} shows 
    \begin{align*}
    N^*(S;X) = \frac{1}{2^{m+n}} \int_{h \in G_0} \#[  S \cap (\Siegel \Lambda h L)_{<X}] \,dh,
    \end{align*}
    with the caveat that elements on the right hand side are weighted by a function similar to \cite[\S6.5, Equation (6.5)]{laga-f4paper}.
    We use the above expression to define $N^*(S;X)$ for subsets $S\subset V(\Z)^{(m,\tau)}$ that are not necessarily $G(\Z)$-invariant.
    Similarly we define $N(S;X) = N^*(S^{\irr};X)$ for every subset $S\subset V(\Z)^{(m,\tau)}$.
    A change of variables trick \cite[Equation (22)]{BhargavaGross} shows that for every $S\subset V(\Z)^{(m,\tau)}$ we then have
    \begin{align*}
    N^*(S;X) = \frac{1}{2^{m+n}} \int_{g\in \Siegel} \int_{\lambda \in \Lambda}\#[S \cap (g\lambda G_0 L)_{<X}] m(g)^{-1} \,dg \, d^{\times}\lambda,
    \end{align*}
    with the caveat that an element $v\in S\cap (g\lambda G_0 L)$ on the right hand side is counted with multiplicity $\#\{h\in G_0 \colon v\in g\lambda h L \}$.

    Recall that $V$ is the space of trace zero linear maps $W\rightarrow W$ that are self-adjoint with respect to $J$.
    Using the basis $e_{-1},\dots,e_{-g},e_0,e_g,\dots,e_1$ of $W$, we may view $V$ as a subspace of matrices. 
    We let $V^{\text{cusp}}$ be the subset of all $T\in V(\R)$ such that the top right entry has absolute value $<1$, i.e. such that $T(e_1) = \sum b_i e_i$ with $|b_{-1}|<1$.
    Let $S_U' = S_U \cap (V(\Z)\setminus  V^{\text{cusp}})$.
    Then
    \begin{align}\label{eq_proofequidist1}
        N(S_U;X) = N(S_U';X) + o(X^{g(2g+3)}) = N^*(S_U';X) + o(X^{g(2g+3)}), 
    \end{align}
    where the first estimate follows from \cite[Proposition 29]{BhargavaGross} (`cutting off the cusp') and the second one follows from \cite[Proposition 31]{BhargavaGross}.
    It remains to estimate $N^*(S_U';X)$.

    We first make the reduction covariant more explicit.
    If $g\in \Siegel, h = (h',\lambda')\in G_0$ (with $h'\in G(\R)$ and $\lambda'\in \Lambda$), $\lambda \in \Lambda$ and $\ell \in L$, the $G(\R)$-equivariance and $\Lambda$-invariance of $\mathcal{R}$ imply that $\mathcal{R}(g h \lambda \ell) = g h' \mathcal{R}(\ell)$.
    We have chosen $L$ so that $\mathcal{R}(\ell) = 1\in G(\R) / K$, so $gh' \mathcal{R}(\ell) = gh'$ in $X_G$.
    Write $\bar{U}$ for the preimage of $U$ under the quotient map $G(\R)\rightarrow G(\Z)\backslash X_G$.
    Then we conclude that $\mathcal{R}(gh\lambda \ell) \in U$ if and only if $gh' \in \bar{U}$.

    Let $Z_1 = \{ (g,h) \in \Siegel \times G_0 \mid gh' \in \bar{U}\}$, where we write $h'$ for the projection of $h \in G_0$ onto $G(\R)$.
    By Lemma \ref{lemma_goodsubsetimpliesregionsemialg}, $Z_1$ is semialgebraic.
    Let $Z_2$ be the graph of the action map $Z_1\times \Lambda\times L\rightarrow V(\R)$ sending $(g,h,\lambda,\ell)$ to $gh\lambda v$.
    Since this map is algebraic, $Z_2$ is semialgebraic.
    So is the projection $Z_3$ of $Z_2$ onto $\Siegel\times V(\R)$.
    Let $Z_4 = \{(g,v,X) \in Z_3\times \R_{>0}\mid \height(v)< X\}$; this is again semialgebraic.
    For every $g\in \Siegel$ and $X\in \R_{>0}$ the set $B(g,X)\coloneqq \{v \in V(\R)\mid (g,v,X)\in Z_4\}$ is semialgebraic and equals $(g\Lambda G_0 L)_{<X} \cap \mathcal{R}^{-1}(U)$.
    We view $B(g,X)$ as a multiset where $v\in B(g,X)$ has multiplicity $\#\{(\lambda,h) \in \Lambda\times G_0 \mid v \in g\lambda h L \}$.
    Then $B(g,X)$ is partitioned into finitely many semialgebraic subsets of constant multiplicity.
    Applying Proposition \ref{prop_countlatticepointsbarroero} to $Z_4$ shows that 
    \begin{align*}
        \#[V(\Z) \cap (B(g,X)\setminus V^{\text{cusp}})] = \Vol(B(g,X) \setminus V^{\text{cusp}}) + E(g,X),
    \end{align*}
    where $E(g,X)$ is the error term. 
    The proof now proceeds in an identical way to that of \cite[Proposition 33]{BhargavaGross}: estimating $E(g,X)$ and $\Vol(B(g,X)\cap V^{\text{cusp}})$ shows that
    \begin{align}\label{eq_proofequidist2}
        N^*(S_U;X) = \frac{1}{2^{m+n}} \int_{g\in \Siegel} \Vol(B(g,X))m(g)^{-1} \, dg + o(X^{g(2g+3)}).
    \end{align}
    The restriction of the characteristic polynomial map $\Lambda \cdot L \rightarrow I(m)$ is a semialgebraic isomorphism; write its inverse by $s\colon I(m)\rightarrow \Lambda\cdot L$.
    Equip $B(\R)$ with the measure corresponding to the differential form $dc_2\wedge dc_3 \wedge \cdots \wedge dc_{2g+1}$. This induces a measure on the open subset $I(m) \subset B(\R)$.
    The change of measure formula of \cite[Proposition 34]{BhargavaGross} shows that 

    \begin{align}\label{eq_proofequidist3}
     \frac{1}{2^{m+n}} \int_{g\in \Siegel} \Vol(B(g,X))m(g)^{-1} \,dg = \frac{|\mathcal{J}|}{2^{m+n}} \Vol(I(m)_{<X})\int_{g \in \Siegel} \int_{h\in G_0}  \mathbf{1}_{\{gh \in \bar{U}\}} m(g)^{-1} \,dh \,dg,
    \end{align}
    where $\mathbf{1}_T$ denotes the indicator function of a set $T$.
    By switching the order of integration, Lemma \ref{lemma_integral_siegelset} and using $\Vol(G_0) = 1$, we calculate that
    \begin{align}\label{eq_proofequidist4}
        \int_{g \in \Siegel} \int_{h\in G_0}  \mathbf{1}_{\{gh \in \bar{U}\}}m(g)^{-1} \, dh \, dg = \int_{h\in G_0} \int_{g\in \Siegel} \mathbf{1}_{\bar{U}h^{-1}}m(g)^{-1} \,  dg \, dh = \int_{h\in G_0} \mu_G(U) \, dh= \mu_G(U).
    \end{align}
    Combining \eqref{eq_proofequidist1}, \eqref{eq_proofequidist2}, \eqref{eq_proofequidist3} and \eqref{eq_proofequidist4} shows that
    \begin{align*}
    N(S_U;X) = \mu_G(U)\frac{|\mathcal{J}|}{2^{m+n}} \Vol(I(m)_{<X}) +o(X^{g(2g+3)})=  \mu_G(U)\frac{|\mathcal{J}|}{2^{m+n}} \Vol(I(m)_{<1})X^{g(2g+3)}+o(X^{g(2g+3)}),
    \end{align*}
    as required.
\end{proof}

\begin{proof}[Proof of Theorem \ref{theorem: equidistribution reduction covariant}]
    Let 
    \[
    \underline{\nu}(U) = \mu_G(G(\Z)\backslash X_G) \liminf_{X\rightarrow \infty} \frac{N(V(\Z) \cap \mathcal{R}^{-1}(U);X)}{N(V(\Z);X)}
    \]
    and $\bar{\nu}(U)$ be the same expression with $\liminf$ replaced by $\limsup$.
    It suffices to prove that $\underline{\nu}(U) = \bar{\nu}(U) = \mu_G(U)$.
    If $U$ is good, this follows from summing \eqref{equation_equidistribution_realcomponent} over all $(m,\tau)$, Equation \eqref{equation_averageorbitslocalversion}, and the fact that $V^s(\R)$ is partitioned into the subsets $V^{(m,\tau)}$.
    To prove the theorem for general $U$, we bootstrap from the case of good subsets (and the trivial case $U = G(\Z) \backslash X_G$).

    Let $U^{\circ}$ be the interior of $U$.
    Since $U^{\circ}$ is open, by Lemma \ref{lemma_countablebasis_goodsubsets} there exists an increasing sequence $(U_n)_{n\geq 1}$ of good subsets whose union is $U^{\circ}$.
    For every $n\geq 1$ we have $\underline{\nu}(U_n) \leq \underline{\nu}(U^{\circ})$.
    Since $U_n$ is good, $\underline{\nu}(U_n) = \mu_G(U_n)$, hence $\mu_G(U_n) \leq \underline{\nu}(U^{\circ})$ for all $n\geq 1$.
    By continuity of the measure, $\mu_G(U_n) \rightarrow \mu_G(U^{\circ})$ as $n\rightarrow \infty$.
    We conclude that $\mu_G(U^{\circ})\leq \underline{\nu}(U^{\circ})$.
    Let $\bar{U}$ denote the closure of $U$.
    Since the complement of $\bar{U}$ equals the interior of the complement of $U$, the above argument also shows $\bar{\nu}(\bar{U}) \leq \mu_G(\bar{U})$.

    In conclusion, we have shown that 
    \[
    \mu_G(U^{\circ}) \leq \underline{\nu}(U^{\circ})\leq \underline{\nu}(U) \leq \bar{\nu}(U) \leq \bar{\nu}(\bar{U}) \leq \mu_G(\bar{U}).
    \]
    By assumption, $\mu_G(U^{\circ}) = \mu_G(\bar{U}) = \mu_G(U)$, so all the inequalities are in fact equalities, and the theorem follows.
\end{proof}

\subsection{A neighbourhood of the cusp}\label{subsec_neighbourhoodcusp}

In our application, we would like to apply Theorem \ref{theorem: equidistribution reduction covariant} with the set $U \subset G(\Z) \backslash X_G$ taken to correspond to the set of of inner products $\langle \cdot, \cdot \rangle_H$ such that $W(\Z)$ contains a short vector, in some sense. In this section, we construct sets which play this role. 

Assume for definiteness that the constant $c > 0$ in the definition of $\Siegel = \omega T_c K$ satisfies $c < 1$. If $\epsilon > 0$, let us define
\[ T(\epsilon) \coloneqq \{ (t_1, \dots, t_g) \in T_{c} \colon t_1 > 1/\epsilon \}, \]
and $U_\epsilon \coloneqq \varphi(\omega T(\epsilon) K) \subset G(\Z) \backslash X_G$, where $\varphi\colon \Siegel\rightarrow G(\Z)\backslash X_G$ denotes the quotient mapping. 
The following proposition contains the essential properties of this set.
\begin{proposition}\label{prop_good_neighbourhoods_of_cusp}
    \begin{enumerate}
        \item The boundary of $U_{\epsilon}$ has measure zero, and $\mu_G(U_\epsilon) \to 0$ as $\epsilon \to 0$.
        \item Suppose that $0 < \epsilon < 1$, that $\langle \cdot, \cdot \rangle_H \in X_G$ and that there exists a nonzero vector $w \in W(\Z)$ such that $\langle w, w \rangle_H^{1/2} < c^g \epsilon$. Then  $G(\Z) \cdot \langle \cdot, \cdot \rangle_H \in U_\epsilon$.
    \end{enumerate}
\end{proposition}
\begin{proof}

    Since the boundary of a semialgebraic set has strictly smaller dimension \cite[Proposition 2.8.13]{BCR-realalgebraicgeometry} and $\varphi^{-1}(U_{\epsilon})$ is semialgebraic, the boundary of $\varphi^{-1}(U_{\varepsilon})$ has measure zero, hence the same is true for $U_{\epsilon}$.
    To show that $\mu_G(U_\epsilon) \to 0$ as $\epsilon \to 0$, it suffices to show that $\Vol( \omega T(\epsilon) K ) \to 0$ as $\epsilon \to 0$. This is true because $\vol(\Siegel) < \infty$ and $\cap_{\epsilon > 0} T(\epsilon) = \emptyset$. 

    To prove the second part, we can assume (since $\varphi$ is surjective) that $\langle \cdot, \cdot \rangle_H = (nt) \cdot \langle \cdot, \cdot \rangle_{H_0}$ for some $n \in \omega$, $t = (t_1, \dots, t_g) \in T_c$. We need to show that the existence of the vector $w$ implies that $t \in T(\epsilon)$, or equivalently that $t_1 >  1/\epsilon$. We at least have 
    \[ t_g > c, t_{g-1} > c^2, \dots, t_1 > c^g, \]
    since $t \in T_c$.
    
    Let $f_{-1}, \dots, f_{-g}, f_0, f_g, \dots, f_1$ denote the image of the standard basis $e_{-1}, \dots, e_1$ under $n t \in P(\R)$. Then the elements $f_{-1}, \dots, f_1$ are orthonormal with respect to $\langle \cdot, \cdot \rangle_H$ and for each $i = 1, \dots, g$ we have
    \[ f_{-i} \in t_{-i} e_{-i} + \langle e_{1-i}, \dots, e_1 \rangle_\R, \]
    \[ f_0 \in e_0 + \langle e_g, \dots, e_1 \rangle_\R, \]
    \[ f_i \in t_i e_i + \langle e_{i+1}, \dots, e_1 \rangle_\R \]
    (where we define $t_{-i} = t_i^{-1}$). We first claim that $w \in \langle e_0, \dots, e_g \rangle_\R$. Indeed, if this is not the case that we can write $w = \sum_{i=j}^g m_{-i} e_{-i} + m_0 e_0 + \sum_{i=1}^g m_i e_i$ for some integers $m_i$ and $j \geq 1$ with $m_{-j} \neq 0$. Then we find that $w \in m_{-j} t_{-j}^{-1} f_{-j} + \langle f_{-(j+1)}, \dots, f_0, f_g, \dots, f_1 \rangle_\R$. Using the orthonormality of the $f_i$, we find that
    \[ \langle w, w \rangle_H^{1/2} \geq |m_{-j}| t_j > c^g, \]
    noting that $|m_{-j}| \geq 1$ since $m_{-j}$ is a nonzero integer. This contradicts the assumption that $\langle w, w \rangle_H^{1/2} < c^g \epsilon < c^g$. A similar argument show that in fact $w \in \langle e_1, \dots, e_g \rangle_\R$. We can therefore write $w = \sum_{i=1}^j m_i e_i$ for some integers $m_i$ and $1 \leq j \leq g$ with $m_j \neq 0$. Then $w \in m_j t_j^{-1} f_j + \langle f_{j-1}, \dots, f_1 \rangle_\R$, and using the orthonormality of the $f_i$ again we find that
    \[ \langle w, w \rangle_H^{1/2} \geq t_j^{-1}. \]
    Since $t_1 > c^g t_j$, we find that 
    \[ c^g \epsilon >  \langle w, w \rangle_H^{1/2} \geq t_j^{-1} > c^g t_1^{-1}, \]
    and therefore that $t_1 > 1/\epsilon$, as required. 
\end{proof}

\begin{corollary}\label{corollary_polyswithsmallnormvectorrare}
    There exist constants $c,c'>0$ such that for every $0< \epsilon <1$,
    \begin{align}\label{equation: lim sup less than mu(U)}
        \limsup_{X\rightarrow \infty} \frac{\#\{f\in \mathcal{F}(X) \colon \exists\, T\in V_f(\Z)^{\irr} \text{ and } w\in W(\Z) \text{ with } \langle w,w\rangle_{H_T}^{1/2} <c^g \epsilon \} }{ \#\mathcal{F}(X)} \leq c'\cdot \mu_G(U_{\epsilon}).
    \end{align}
    Moreover, $\mu_G(U_{\epsilon})\rightarrow 0$ as $\epsilon \rightarrow 0$.
\end{corollary}
\begin{proof}
    Denote the numerator in the above $\limsup$ by $\mathcal{F}(X)^{\text{short}}$.
    Part 2 of Proposition \ref{prop_good_neighbourhoods_of_cusp} shows that $
    \#\mathcal{F}(X)^{\text{short}} \leq N(V(\Z)\cap \mathcal{R}^{-1}(U_{\epsilon});X)$.
    Theorem \ref{theorem: equidistribution reduction covariant} shows that 
    \[ N(V(\Z)\cap \mathcal{R}^{-1}(U_{\epsilon});X) =  \frac{\mu_G(U_{\epsilon})}{\mu_G(G(\Z)\backslash X_G)} N(V(\Z);X) + o(X^{g(2g+3)}). \]
    We also have $N(V(\Z);X) = c_1 \#\mathcal{F}(X) + o(X^{g(2g+3)})$ for some $c_1>0$, by \eqref{equation_averageirredorbitsfinite}. Setting $c' =\mu_G(G(\Z)\backslash X_G)^{-1} c_1$, we see that $\#\mathcal{F}(X)^{\text{short}} \leq c'\cdot \mu_G(U_{\epsilon})  \#\mathcal{F}(X) + o(X^{g(2g+3)})$, as desired.
\end{proof}

\section{Classification and construction of orbits}\label{sec_classificationconstructionorbits}

In this section we first review in \S\ref{subsec_rationalorbits} the construction of elements of $G(\Q) \backslash V(\Q)$ from rational points of $J_f(\Q)$, following \cite{BhargavaGross} and \cite{Thorne-remark}, before giving our constructions of refined integral structures and discussing how they interact with the reduction covariant in \S\ref{subsec_integralorbits}.
We then study in \S\ref{subsec_canonicalplotdistinguishedorbit} the canonical plot of the reduction covariant of the integral orbit associated to the identity in the Jacobian.
Finally in \S\ref{subsec_Qinvariant} we study more generally the reduction covariants of integral orbits that lie in the distinguished $G(\Q)$-orbit. This will be useful in the study of the reduction covariant associated to points in $2J_f(\Q)$.

\subsection{Rational orbits}\label{subsec_rationalorbits}

Let $k / \Q$ be a field, and let $W_0$ be a quadratic space over $k$ of rank $2g+1$, as in \S\ref{subsec_quadspaces}.  We say that $W_0$ is \define{split} if there exists a $g$-dimensional subspace $L_0 \leq W_0$ such that the restriction of $(\cdot, \cdot )_{W_0}$ to $L_0 \times L_0$ is 0. An isomorphism $W_0 \to W_0'$ of quadratic spaces is, by definition, a $k$-linear map intertwining the two symmetric bilinear forms.
\begin{lemma}\label{lem_split_spaces_isomorphic}
Let $W_0, W_0'$ be split quadratic spaces of the same discriminant. Then they are isomorphic (as quadratic spaces over $k$).
\end{lemma}
\begin{proof}
If $W_0$ is split, then there is an isomorphism of quadratic spaces $W_0 \cong H^g \oplus \langle a \rangle$, where $H$ is a hyperbolic plane and $\langle a \rangle$ is the 1-dimensional quadratic space which represents $a$, for some $a \in k^\times$. We now observe that the value of $a$ is determined by $\disc W_0$.
\end{proof}

\begin{proposition}\label{prop_classification_of_split_orbits}
Let $f(x) = x^{2g+1} + c_2 x^{2g-1} + \dots + c_{2g+1} \in k[x]$ be a separable polynomial. Then the map $T \mapsto (W, T)$ determines a bijection between the following two sets:
\begin{enumerate} \item The set of $G(k)$-orbits in $V_f(k)$.
\item The set of isomorphism classes of pairs $(W_0, T_0)$, where $W_0$ is a split quadratic space over $k$ of discriminant $1 \in k^\times / (k^\times)^2$ and $T_0 : W_0 \to W_0$ is a self-adjoint linear operator of characteristic polynomial $f(x)$.
\end{enumerate}
\end{proposition}
\begin{proof}
    Use Lemma \ref{lem_split_spaces_isomorphic} and the fact that $\mathrm{O}(W_0) = \SO(W_0) \times \{ \pm 1\}$ since $\dim_k W_0$ is odd; see \cite[Lemma 2.2]{Thorne-remark}.
\end{proof}
This proposition is useful in conjunction with the following construction. If $f(x) = x^{2g+1} + c_2 x^{2g-1} + \dots + c_{2g+1} \in k[x]$, then we define $A_f = k[x] / (f(x))$, and write $\tau : A_f \to k$ for the linear form which sends a polynomial of degree $\leq 2g$ to the  coefficient of $x^{2g}$. If $U \in A_f^\times$, then we can define a quadratic space $W_U$ and self-adjoint linear operator $T_U : W_U \to W_U$ by defining $(v, w)_{W_U} = \tau(vw / U)$ and $T_U(v) = x v$.
\begin{proposition}\label{prop_classification_of_orbits}
Let $f(x) = x^{2g+1} + c_2 x^{2g-1} + \dots + c_{2g+1} \in k[x]$ be a separable polynomial. Then the map $U \mapsto (W_U, T_U)$ determines a bijection between the following two sets :
\begin{enumerate}
    \item The group $A_f^\times / (A_f^\times)^2$.
    \item The set of isomorphism classes of pairs $(W_0, T_0)$, where $W_0$ is a quadratic space over $k$ and $T_0 : W_0 \to W_0$ is a self-adjoint linear operator of characteristic polynomial $f(x)$.
\end{enumerate}
Under this bijection, a class $U \text{ mod }(A_f^\times)^2$ corresponds to a quadratic space of discriminant $\mathbf{N}_{A_f / k}(U) \in k^\times / (k^\times)^2$. 
\end{proposition}
\begin{proof}
See \cite[Proposition 2.5]{Thorne-remark}. 
\end{proof}
Combining the above two propositions, we see that we can construct rational orbits in $G(k) \backslash V_f(k)$ by finding elements $U \in A_f^\times$ of square norm such that the associated quadratic space $W_U$ is split. One such element is the class $1 \in A_f^\times$ of the identity. Indeed, in this case $L_1 = \langle 1, x, \dots, x^{g-1} \rangle \leq W_1$ is an isotropic subspace of dimension $g$. We see moreover that $L_1^\perp = \langle 1, x, \dots, x^g \rangle$ and therefore that $T_1 L_1 \leq L_1^\perp$. We have therefore constructed the distinguished orbit of $V_f(k)$. 

Following \cite{Thorne-remark}, we next show how we can find other  elements $U$ such that $W_U$ is split using the Mumford representation of divisors on hyperelliptic curves. Take a separable polynomial $f(x) = x^{2g+1} + c_2 x^{2g-1} + \dots + c_{2g+1} \in k[x]$. Then the affine curve
\[ C_f^0 : y^2  = f(x) \]
is smooth, and has a smooth projective compactification $C_f$, which is a hyperelliptic curve of genus $g$ with a unique marked point $P_\infty \in C_f(k) \setminus C_f^0(k)$. We set $J_f = \Pic^0_{C_f}$. Any divisor class $[D] \in J_f(k)$ has a unique representative of the form $D = \sum_{i=1}^m P_i - m P_\infty$, where the points $P_i \in C_f(k^s)$ (defined over a separable closure $k^s / k$ and unique up to re-ordering), satisfy:
\begin{itemize}
    \item $0 \leq m \leq g$; 
    \item $P_i \neq P_\infty$ for each $i = 1, \dots, m$;
    \item $P_i \neq \iota(P_j)$ for each $i \neq j$, where $\iota : C_f \to C_f$ denotes the hyperelliptic involution (defined by $\iota(y) = -y$);
\end{itemize}
see \cite[Ch. IIIa, \S1]{MumfordTataII}. We call $m$ the \define{Mumford degree} of the divisor class $[D]$. Associated to the divisor $D$ is its Mumford representation $(U, V, R)$: a tuple of polynomials in $k[x]$ satisfying the following conditions:
\begin{itemize}
\item $U(x) = \prod_{i=1}^m(x - x(P_i))$. In particular, $U$ is monic of degree $m$.
\item $R(x)$ has degree $\leq m-1$ and satisfies $R(x(P_i)) = y(P_i)$ for each $i = 1, \dots, m$ (interpreted with multiplicity if $x(P_i)$ is a repeated root of $U$).
\item $V(x)$ is monic of degree $2g+1-m$ and $UV = f - R^2$. 
\end{itemize} 
The assignment $[D] \mapsto (U, V, R)$ determines, for each $m = 0, \dots, g$, a bijection between the set of points $[D] \in J_f(k)$ of Mumford degree $m$ and the set of tuples $(U, V, R)$ of polynomials in $k[x]$ such that $U$ is monic of degree $m$, $V$ is monic of degree $2g + 1 - m$, $R$ has degree $\leq m-1$, and $UV = f - R^2$. We call a tuple $(U, V, R)$ satisfying these conditions a \define{Mumford triple} for $f$.

The space $H^0( C_f^0, \cO_{C_f}(-D))$ is a locally free $k[x, y] / (y^2 - f(x))$-module of rank 1. Therefore $W_D = H^0( C_f^0, \cO_{C_f}(-D)) / (y)$ is a free $A_f$-module of rank 1. We can give it the structure of quadratic space as follows: there is an isomorphism 
\[ \cO_{C_f}(-D) \otimes_{\cO_{C_f}} \cO_{C_f}(-\iota^\ast D) \cong \cO_{C_f}( - (D + \iota^\ast D) ) \cong \cO_{C_f}, \]
uniquely characterised by the requirement that the nonvanishing global section $U \in H^0(C_f, \cO_{C_f}( - (D + \iota^\ast D) ) )$ is sent to $1 \in H^0(C_f, \cO_{C_f})$. After quotienting by $y$, this determines an isomorphism $\varphi_D : W_D \otimes_{A_f} W_D \cong A_f$, and we define $(v, w)_{W_D} = (-1)^{\deg U} \tau( \varphi_D(v \otimes w) )$. Observe that $W_D$ comes equipped with a self-adjoint linear endomorphism $T_D$ of characteristic polynomial $f(x)$, given by the formula $T_D(v) = xv$.
\begin{proposition}
Let $f(x) = x^{2g+1} + c_2 x^{2g-1} + \dots + c_{2g+1} \in k[x]$ be a separable polynomial, and let $(U, V, R)$ be a Mumford triple for $f$. Factorise $U(x) = U_0(x) U_1(x)$, where $U_0(x), U_1(x) \in k[x]$ are monic, $U_0 | f$, and $(U_1, f) = 1$. Then under the bijection of Proposition \ref{prop_classification_of_orbits}, the pair $(W_D, T_D)$ corresponds to the class $(-1)^{\deg U} U_1 (U_0 - f / U_0) \text{ mod }(A_f^\times)^2$. 

Moreover, we have $\disc W_D = 1 \text{ mod }(k^\times)^2$ and the quadratic space $W_D$ is split, so $(W_D, T_D)$ determines (via the bijection of Proposition \ref{prop_classification_of_split_orbits}) an element of $G(k) \backslash V_f(k)$. 
\end{proposition}
We remark that the generic case occurs when $U_0 = 1$, hence $U_1 (U_0 - f / U_0) = U$ in $A_f$. 
\begin{proof}
We review the proof from \cite{Thorne-remark}. The space $H^0( C_f^0, \cO_{C_f}(-D))$ is a free $k[x]$-module with basis $U, y-R$. Using the relations $U V = f - R^2$ and $y^2 = f$, hence $(y-R)(y+R) = U V$, we see that multiplication by $y$ acts by the $2 \times 2$ matrix
\[ \left( \begin{array}{cc} R& V \\ U & -R \end{array} \right). \]
Writing $u = \deg U$, $v = \deg V$, the elements $U, x U, \dots x^{v-1} U$, $(y-R), x(y-R), \dots, x^{u-1} (y-R)$ project to a $k$-basis of $W_D$. We have $\varphi_D( U \otimes U) = U$, $\varphi_D (U \otimes (y-R)) = -R$, $\varphi_D( (y-R) \otimes (y-R) ) = -V$.

Let $a = \lfloor u / 2 \rfloor$, $b = \lfloor v /2  \rfloor$. Using the above relations, one can compute that the elements $U, \dots, x^{b-1} U$, $(y-R), x(y-R), \dots, x^{a-1} (y-R)$ form a basis for a $g$-dimensional isotropic subspace $L_D$ of $W_D$. In particular, $W_D$ is split. If $u$ is even then $L_D^\perp = \langle L_D, x^b U \rangle_k$ and $\disc W_D = (x^b U, x^b U)_{W_D} = \tau( x^{v-1} U ) = 1$. If $u$ is odd then $L_D^\perp = \langle L_D, x^a (y-R) \rangle_k$ and $\disc W_D = (x^a (y-R), x^a (y-R))_{W_D} = \tau( x^{u-1} V) = 1$. 

To complete the proof, we need to show that $(W_D, T_D)$ corresponds to the class $(-1)^{\deg U}U_1 (U_0 - f / U_0) \text{ mod }(A_f^\times)^2$. This class is computed as follows: choose a vector $w \in W_D$ which generates $W_D$ as an $A_f$-module (where $x \in A_f$ acts as $T_D$), and compute $(-1)^{\deg U}\varphi_D(w \otimes w)$. We can take $w = U_1( y - U_0 ) $, giving 
\[ \varphi_D(w \otimes w) = U_1^2 (y-U_0)(-y-U_0) / U  = U_1^2 (U_0^2 - f) / U_0 U_1  = U_1 (U_0 - f / U_0), \]
as claimed. 
\end{proof}
\begin{corollary}\label{corollary_rationalorbits_exist}
    Let $f(x) = x^{2g+1} + c_2 x^{2g-1} + \dots + c_{2g+1} \in k[x]$ be a separable polynomial. Then the map $[D] \mapsto (W_D, T_D)$ determines an injection $J_f(k) / 2 J_f(k) \to G(k) \backslash V_f(k)$, which sends the identity to the distinguished orbit. 
\end{corollary}
\begin{proof}
The composite $J_f(k) \to G(k) \backslash V_f(k) \to A_f^\times / (A_f^\times)^2$ sends $[D]$ to $(-1)^{\deg U} U_1 (U_0 - f / U_0) \text{ mod }(A_f^\times)^2$. It follows from \cite[Lemma 2.2]{Thorne-remark} that this map agrees with the Kummer 2-descent map, and in particular that if factors through an injection with source $J_f(k) / 2 J_f(k)$. (We remark that \cite{Thorne-reduction} takes the descent map as having image in $A_f^\times / k^\times (A_f^\times)^2$, whereas here we work with a map taking image in $\ker( \mathbf{N}_{A_f / k} : A_f^\times / (A_f^\times)^2 \to k^\times / (k^\times)^2 )$. The cited result  gives what we need because the map 
\[ \ker( \mathbf{N}_{A_f / k} : A_f^\times / (A_f^\times)^2 \to k^\times / (k^\times)^2 ) \to A_f^\times / k^\times (A_f^\times)^2 \]
is an isomorphism, because $[A_f : k] = 2g+1$ is odd.)
The identity of $J_f(k)$ corresponds to the Mumford triple $(1, f, 0)$, which indeed is mapped to the identity of $A_f^\times / (A_f^\times)^2$, which corresponds to the distinguished orbit. 
\end{proof}
\subsection{Integral orbits}\label{subsec_integralorbits}

Let $R$ be a Dedekind domain of characteristic 0, and let $k = \Frac(R)$. Let $f = x^{2g+1} + c_2 x^{2g-1} + \dots + c_{2g+1} \in R[x]$ be a polynomial of nonzero discriminant. Then there is a natural map
\[ G(R) \backslash V_f(R) \to G(k) \backslash V_f(k) \]
which is, in general, neither injective nor surjective. Nevertheless, a natural way to construct elements of $G(R) \backslash V_f(R)$ is by constructing $k$-orbits and then trying to give them an $R$-structure. We call an \define{$R$-lattice} in a quadratic space $W_0$ over $k$ an $R$-submodule $M_0 \leq W_0$, generated by a $k$-basis for $W_0$, such that $( \cdot, \cdot)_{W_0}|_{M_0 \times M_0}$ takes values in $R$ and induces an isomorphism $M_0 \cong \Hom_R(M_0, R)$.
\begin{proposition}\label{prop_classification_of_integral_lattices_in_split_quadratic_spaces}
    Suppose that $R = \Z_p$ for some prime $p$, or that $R = \Z$, and let $k = \Frac R$. Suppose that $W_0$ is split, and that $M_0 \leq W_0$ is an $R$-lattice. Then there exists an isomorphism $W_0 \to W(k)$ of quadratic spaces that restricts to an isomorphism $M_0 \cong W(R)$.
\end{proposition}
\begin{proof}
If $R = \Z_p$, this is the content of \cite[Lemmas 15, 20]{BhargavaGross}. If $R = \Z$, then it follows from \cite[Ch. V, Theorem 4]{SerreCiA}. 
\end{proof}

Before stating the classification of integral orbits in $V$, we prove the following proposition, which will be used in \S\ref{subsec_Qinvariant}.

\begin{proposition}\label{prop_uniqueness_of_cusp}
    Suppose that $R = \Z_p$ for some prime $p$, or that $R = \Z$, and let $k = \Frac R$. Then $G(R) \cdot P(k) = G(k)$. Consequently, for any flag of $k$-vector subspaces
    \[ 0 = F_0 \leq F_1 \leq F_2 \dots \leq F_{2g+1} = W(k)\]
    such that $\dim_k F_i = i$ and $F_i^\perp = F_{2g+1-i}$ for each $i = 0, \dots, 2g+1$, we can find an $R$-basis \[ b_{-1}, \dots, b_{-g}, b_0, b_g, \dots, b_1 \]
    for $W(R)$ satisfying
    \[ F_i = \langle b_1, \dots, b_i \rangle_k\text{ } (1 \leq i \leq g), \]
    \[ F_{g+1} = \langle b_0, \dots, b_g \rangle_k, \]
    \[ F_{2g+2-i} = \langle b_{-i}, \dots, b_{-g}, b_0, \dots, b_g \rangle_k \text{ }(1 \leq i \leq g), \]
    and $(b_i, b_j)_W = \delta_{i, -j}$ for each $-g \leq i, j \leq g$.
\end{proposition}
\begin{proof}
    We just give the proof in the harder case that $R = \Z$. In this case the given statement is equivalent to the assertion that $G(\Z)$ acts transitively on the set of flags $F_\bullet$ in $W(\Q)$ as in the statement of the proposition.

    We first choose any $R$-basis $f_{-1}, \dots, f_{-g}, f_0, f_g, \dots, f_1$ of $W(R)$ such that $F_i$ is spanned by the last $i$ basis vectors in the list for each $i = 1, \dots, 2g+1$. Let $A$ denote the Gram matrix of $(\cdot,\cdot)_W$ with respect to this choice of basis. Then the entries of $A$ below the anti-diagonal are 0, the other entries lie in $R$, and $\det A = \det J = (-1)^g$ (noting that any unit in $\Z$ has square $1$). After possibly multiplying some of the $f_i$'s by $\pm 1$, we can assume that the anti-diagonal entries of $A$ are all equal to $1$. We need to show that we can find a unipotent lower-triangular matrix $n \in \GL_{2g+1}(\Z)$ such that ${}^t n A n = J$. Indeed, this would imply the existence of an element $\gamma \in G(\Z)$ taking the standard flag $F_{0, \bullet}$ of $W$ (where $F_{0, i}$ is spanned by the last $i$ basis vectors in the list $e_{-1}, \dots, e_{-g}, e_0, e_g, \dots, e_1$) to $F_\bullet$. 
    
     It is easy to show that we can find an $n$ such that 
    \[ {}^t n A n =  \left( \begin{array}{ccc} C \Psi + B + \Psi {}^t C + y {}^t y & 0 & \Psi \\ 0 & 1 & 0 \\ \Psi & 0 & 0 \end{array} \right), \]
    where $\Psi$ is the $g \times g$ matrix with $1$'s on the anti-diagonal and $0$'s elsewhere, $B \in \Mat_g(\Z)$ is a symmetric matrix which depends on $A$, and $C \in \Mat_g(\Z)$, $y \in \Z^g$ may be chosen freely. We need to show that, given any $B$, we can choose $C, y$ so that $C \Psi + B + \Psi {}^t C + y {}^t y = 0$. We can take $y = (b_{11}, \dots, b_{gg})$, so that $y {}^t y = ( y_i y_j )$ and $B + y {}^t y$ is a symmetric matrix with even diagonal entries; then we can choose $C$ with the required property. 
\end{proof}
\begin{proposition}\label{prop_classification_integralorbits}
    Let $R = \Z_p$ for some prime $p$, or $R = \Z$, and let $f = x^{2g+1} + c_2 x^{2g-1} + \dots + c_{2g+1} \in R[x]$ be a polynomial of nonzero discriminant. Then the assignment $T\mapsto (W(k),W(R),T)$ induces a bijection between:
    \begin{enumerate}
    \item The set of $G(R)$-orbits in $V_f(R)$.
    \item The set of isomorphism classes of tuples $(W_0, M_0, T_0)$, where $W_0$ is a split quadratic space over $k$ of discriminant $1 \text{ mod }(k^\times)^2$, $T_0$ is a self-adjoint $k$-linear endomorphism of $W_0$ of characteristic polynomial $f(x)$, and $M_0 \leq W_0$ is an $R$-lattice such that $T_0 M_0 \leq M_0$.
    \end{enumerate}
\end{proposition}
\begin{proof}
We construct an inverse. Given a tuple $(W_0, M_0, T_0)$, Proposition \ref{prop_classification_of_integral_lattices_in_split_quadratic_spaces} shows that we can find an isomorphism $W_0 \to W(k)$ taking $M_0$ to $W(R)$, and we obtain a vector in $V_f(R)$ by transporting $T_0$ to an endomorphism $T$ of $W(R)$. We need to check that this is well-defined. Any other isomorphism $W_0 \to W(k)$ taking $M_0$ to $W(R)$ differs from the given one by an element of $\mathrm{O}(W)(R) = G(R) \times \{ \pm 1 \}$, so the $G(R)$-orbit of $T$ is indeed independent of any choices made. 
\end{proof}
If $W_0$ is a quadratic space over $\Q$, then giving a $\Z$-lattice $M_0 \leq W_0$ is equivalent to giving a $\Z_p$-lattice $M_{0, p}$ in $W_0 \otimes_\Q \Q_p$ for every prime $p$ (subject to the usual boundedness constraint). This means that an effective way to construct elements of $G(\Z) \backslash V_f(\Z)$ is to construct a pair $(W_0, T_0)$ over $\Q$ and then specify an integral structure prime-by-prime. The main input in such a construction is the following key lemma.
\begin{lemma}\label{lem_existence_of_good_local_lattice}
    Let $p$ be a prime, and let $f = x^{2g+1} + c_2 x^{2g-1} + \dots + c_{2g+1} \in \Z_p[x]$ be a polynomial of nonzero discriminant. Let $(U, V, R)$ be a Mumford triple for $f$ over $\Q_p$, and let $r \geq 0$ be minimal such that $p^{2r} U \in \Z_p[x]$. Then there exists a $\Z_p$-lattice $M_D \leq W_D$ such that $T_D M_D \leq M_D$ and $p^r U \in  M_D - p M_D$. 
\end{lemma}

Note that $p=2$ is allowed in Lemma \ref{lem_existence_of_good_local_lattice}.
This shows that the divisibility condition in \cite[Proposition 24]{BhargavaGross} is unnecessary.

\begin{proof}
    Let $R_f = \Z_p[x] / (f(x)) \leq A_f$. We first say a bit more about $r$. We can factor $U(x) = U_+(x) U_-(x)$, where $U_+ \in \Z_p[x]$ and all the slopes of the Newton polygon of $U_-$ are (strictly) negative. Then the constant term of $U_-$ is the unique coefficient of smallest valuation. If $k \geq 0$ is minimal such that $p^k U_- \in \Z_p[x]$, then we have $p^k U_- \text{ mod }(f) \in R_f^\times$ (as the reduction modulo $p$ of this element lies in the image of $\F_p^\times \to (R_f / (p))^\times$, hence is a unit). Moreover, $k$ is even, because the norm of $p^k U_- \text{ mod }(f)$ is a $p$-adic unit, while the norm of $U_- \text{ mod }(f)$ is a square in $\Q_p^\times$. Thus $k = 2r$ and $p^{2r} U_- \in \Z_p[x]$ has the property that its image in $R_f$ is a unit. 

    We now show that if the lemma holds in the case $r = 0$, then it holds in the general case. Suppose that $r > 0$, hence $U = U_+ U_-$ with $\deg U_- > 0$. We can correspondingly write $D = \sum_i P_{+, i} + \sum_j P_{-, j} - m P_\infty$ for points $P_{\pm, i} \in C_f^0(\overline{\Q}_p)$, and we let $D_+ = \sum_i P_{+, i} - \deg U_+ \cdot P_\infty$. Then there is a canonical isomorphism $W_D \cong W_{D^+}$ of $\Q_p$-vector spaces (because $U_-$ is prime to $f$, so the morphism  of coherent sheaves $\cO_{C_f^0}(-D) \to \cO_{C_f^0}(-D^+)$ is an isomorphism in a neighbourhood of the locus $y = 0$), which identifies $(\cdot, \cdot)_{W_{D+}}$ with the form on $W_D$ given by the formula $(v, w)_+ = (-1)^{\deg U_-} ( U_- v, w )_{W_D}$. 

    By assumption, there exists a $T_{D^+}$-invariant $\Z_p$-lattice $M_{D^+} \leq W_{D+}$ such that $U_+ \in M_{D+} - p M_{D+}$. Pulling back to $W_D$ gives a $T_D$-invariant $\Z_p$-submodule $M_+ \leq W_D$ which is a $\Z_p$-lattice for the form $(\cdot, \cdot)_+$ and such that $U_+ \in M_+ - p M_+$. It follows that $p^{-r} M_+$ is a $T_D$-invariant $\Z_p$-lattice for the form $p^{2r} (\cdot, \cdot)_+ = (-1)^{\deg U_-} ( p^{2r} U_- \cdot, \cdot )_{W_D}$. Since $p^{2r} U_-$ is a unit in $R_f^\times$, hence acts invertibly on $p^{-r} M_+$, it follows that $M_D = p^{-r} M_+$ is a $T_D$-invariant $\Z_p$-lattice for the form $(\cdot, \cdot)_{W_D}$. Moreover, by construction $p^{-r} U_+ \in M_D - p M_D$, hence $p^{-r} U_+ \cdot p^{2r} U_- = p^r U \in M_D - p M_D$, which is the required property. 

    For the rest of the proof, we can therefore assume that $r = 0$. We make some further remarks. If $(U, V, R)$ is a Mumford triple for $f$ where $U \in \Z_p[x]$ has degree $m$ and $R \not\in \Z_p[x]$, then consideration of the relation $UV = f - R^2$, together with the Newton polygon of $f - R^2$, shows that if $s > 0$ is the least integer such that $p^s R \in \Z_p[x]$, and if we factor $V = V_+ V_-$ where $V_+ \in \Z_p[x]$ and the Newton polygon of $V_-$ has (strictly) negative slopes, then $p^{2s} V_-   \in \Z_p[x]$ and $p^{2s} V_- \text{ mod }(f) \in R_f^\times$. Moreover, $\deg V_+ \leq m-2$. 

    To establish the lemma in the case $r = 0$, we will prove a slightly stronger claim, namely:
    \begin{itemize} \item[$(\star)$] If $(U, V, R)$ is a Mumford triple for $f$ with $U \in \Z_p[x]$ of degree $m$, and $s$ is defined as in the previous paragraph, then there exists a $T_D$-invariant $\Z_p$-lattice $M_D \leq W_D$ such that $U, p^s(y-R) \in M_D$.
    \end{itemize} 
    Note that $U$ is then automatically saturated in $M_D$ (i.e. $U \in M_D - p M_D$) as we have $(U, x^{v-1} U)_{W_D} = (-1)^{\deg U}$, which is not possible if $U \in p M_D$.
    
    We first show that $(\star)$ holds when $s = 0$. Indeed, we claim that in this case, we can take $M_D = R_f \cdot U + R_f \cdot (y-R)$. Let $u = \deg U$, $v = \deg V$, and let $M_D'$ denote the $\Z_p$-span of the elements $U, \dots, x^{v-1} U$, $(y-R), \dots, x^{u-1} (y-R)$. Then $M_D' \leq M_D$. We observe that the form $(\cdot, \cdot)_{W_D}$ is $\Z_p$-valued on $M_D$, since we have the formulae
    \[ (x^i U, x^j U)_{W_D} = (-1)^u \tau(x^{i+j} U), \]
    \[ (x^i U, x^j(y-R))_{W_D} = (-1)^{u-1} \tau(x^{i+j} R), \] 
    \[ (x^i (y-R), x^j(y-R))_{W_D} = (-1)^u \tau(x^{i+j} V) \]
    for any $i, j \geq 0$, and we are assuming that $U, V, R \in R_f$. We will show that we can order the basis elements of $M_D'$ so that the Gram matrix of $(\cdot, \cdot)_{W_D}$ has 0's above the anti-diagonal and $\pm 1's$ on the anti-diagonal. Since $M_D' \leq M_D$, this will show at the same time that $M_D'$ is a $\Z_p$-lattice and that $M_D' = M_D$, hence that $M_D$ is a $\Z_p$-lattice which is also $T_D$-invariant. 

    We now describe the basis. Let $a = \lfloor u /2 \rfloor$, $b = \lfloor v / 2 \rfloor$, so $a \leq b$. We take the basis  defined as follows: the elements in each of the two lists $U, \dots, x^{v-1} U$ and $(y-R), \dots, x^{u-1} (y-R)$ will appear in order, but we need to interleave these two lists. We take $(b-a+1)$ elements of the form $x^i U$, then alternately take elements of the form $x^j (y-R)$ and $x^i U$ until the $x^j (y-R)$'s are exhausted, then take the remaining elements of the form $x^i U$. We leave the elementary check that this basis has the claimed properties to the reader. 

    We now prove $(\star)$ in general by induction on $m$, the case $m = 0$ being trivial. Suppose that $0 < m \leq g$ and that $(U, V, R)$ is a Mumford triple for $f$ with $U \in \Z_p[x]$ of degree $m$, $p^s R$ saturated in $\Z_p[x]$ for some $s > 0$. Then we can factor $V = V_+ V_-$ as above. Let $\alpha_i \in \overline{\Q}_p$ be the roots of $V_+$, with multiplicity, and let $\beta_i = -R(\alpha_i)$, $m' = \deg V_+$, $D' = \sum_i (\alpha_i, \beta_i) - m' P_\infty$. The divisor of the function $(y+R)$ is $\iota^\ast D + E$, where $E = D' + D''$ for a divisor $D''$ whose intersection with $C_f^0$ is supported on the zeroes of $V_-$. The divisor of the function $U$ is $D + \iota^\ast D$. Therefore the divisor of $(y+R) / U$ is $E - D$, and multiplication by $(y+R) / U$ defines an isomorphism $W_D \to W_{D'}$ of $\Q_p$-vector spaces which identifies $(\cdot, \cdot)_{W_{D'}}$ with the form $(\cdot, \cdot)'$ on $W_D$ defined by $(v, w)' = (-1)^{\deg V_- - 1}(V_- v, w)_{W_D}$. Let $R' \in \Q_p[x]$ be the unique polynomial of degree $\leq m'-1$ such that $R'(\alpha_i) = \beta_i$ for each root $\alpha_i$ of $V^+$; then we can write $R + R' = V_+ R_1$ for some $R_1 \in \Q_p[x]$. Let $s' \geq 0$ be the least integer such that $p^{s'} R' \in \Z_p[x]$. Then $s' \leq s$, since we can construct $p^s R'$ by subtracting a $\Z_p[x]$-multiple of $V_+$ from $-p^s R$. For the same reason, we have $p^s R_1 \in \Z_p[x]$. 
    
    Applying the induction hypothesis to $D'$ and pulling back to $W_D$, we see that there is a $T_D$-invariant $\Z_p$-submodule $M' \leq W_D$ which is a $\Z_p$-lattice for the form $( \cdot, \cdot )'$ and such that the elements $U V_+ / (y+R)$, $p^{s'} U (y-R') / (y+R)$ are contained in $M'$. It follows that $p^{-s} M'$ is a $T_D$-invariant $\Z_p$-lattice for the form $p^{2 s}( \cdot, \cdot)' = (-1)^{\deg V_- - 1}(p^{2s} V_- \cdot, \cdot)_{W_D}$. Since $p^{2 s}  V_- \in R_f^\times$, it follows that $M_D = p^{-s} M'$ is a $T_D$-invariant $\Z_p$-lattice for the form $(\cdot, \cdot)_{W_D}$. 
    
    To complete the induction step, we need to show that the elements $U, p^s (y-R)$ lie in $M_D$. By construction, $p^{-s} U V_+ / (y+R)$, $p^{s' - s} U (y-R') / (y+R)$ (or rather, the images of these elements of $H^0(C_f^0, \cO_{C_f}(-D))$ in $W_D$) lie in $M_D$. The identity $UV = y^2 - R^2$ gives $U V_+ / (y+R) = (y-R) / V_+$, hence $p^{-s} U V_+ / (y+R) = p^{-s} (y-R) / V_-$. Since $p^{2s} V_- \in R_f^\times$, we find that $p^{2s} V_- \cdot p^{-s} (y-R) / V_- = p^s (y-R) \in M_D$. Similarly, we have $R = -R' + V_+ R_1$, hence $y+R = y-R' + V_+ R_1$, hence $U (y-R') / (y+R) = U - U V_+ R_1 / (y+R)$. Since $s' \leq s$, we have $U (y-R') / (y+R) = U - U  V_+ R_1 / (y+R) \in M_D$. To show that $U \in M_D$, it is therefore enough to show that $U V_+ R_1 / (y+R) \in M_D$. We can write $U V_+ R_1 / (y+R) = (p^s R_1) \cdot (p^{-s} U V_+ / (y+R))$, so this follows since $p^{-s} U V_+ / (y+R) \in M_D$ and $p^s R_1 \in \Z_p[x]$. This completes the proof. 
\end{proof}
\begin{theorem}\label{thm_good_global_lattice}
    Let $f = x^{2g+1} + c_2 x^{2g-1} + \dots + c_{2g+1} \in \Z[x]$ be a polynomial of nonzero discriminant. Let $(U, V, R)$ be a Mumford triple for $f$ over $\Q$ associated to the divisor $D$, and let $M \in \mathbb{N}$ be the least natural number such that $ M U \in \Z[x]$. Then $M = N^2$ is a square. Suppose that $(U, f) = 1$. Then there exists a $\Z$-lattice $M_D \leq W_D$ such that $T_D M_D \leq M_D$ and there is a primitive vector $v_D \in M_D$ satisfying 
    \begin{align}\label{equation_normofvD}
    \frac{1}{2} \log \langle v_D, v_D \rangle_{H_{T_D}} = \frac{1}{2} \log M + \frac{1}{2} \log \sum_{i=1}^{2g+1} \frac{ | U(\omega_i) |}{ | f'(\omega_i) |},
    \end{align}
    where $\omega_1, \dots, \omega_{2g+1} \in \bC$ are the roots of $f(x)$.
\end{theorem}
\begin{proof}
    The first paragraph of Lemma \ref{lem_existence_of_good_local_lattice} shows that $M$ is a square.
    Factor $N = \prod_p p^{r_p}$. By Lemma \ref{lem_existence_of_good_local_lattice}, we can find, for each prime $p$, a $T_D$-invariant $\Z_p$-lattice $M_{D, p} \leq W_D \otimes_{\Q} \Q_p$ such that $p^{r_p} U \in M_{D, p} - p M_{D, p}$. There exists a unique $\Z$-lattice $M_D \leq W_D$ such that for each prime $p$, we have $M_D \otimes_\Z \Z_p = M_{D, p}$. This lattice is necessarily $T_D$-invariant, and $N U \in M_D$ is a primitive vector; we take $v_D = N U$. 

    To prove the theorem, we need to compute the norm of this vector with respect to the inner product $\langle \cdot, \cdot \rangle_{H_{T_D}}$, using the description following Proposition \ref{prop_existence_of_reduction_covariant}. Since $(U, f) = 1$, the inclusion of sheaves $\cO_{C_f^0}(-D) \to \cO_{C_f^0}$ is an isomorphism in a neighbourhood of $y = 0$, and we can identify $W_D = W_1 = A_f$. Consequently, there is an isomorphism $W_D \otimes_\Q \bC \cong A_f \otimes_\Q \bC \cong \bC^{2g+1}$ (where the last isomorphism sends a polynomial $g(x)$ to $(g(\omega_i))_{i=1, \dots, 2g+1}$). Let $e_i \in W_D \otimes_\Q \bC$ denote the pre-image of the $i^\text{th}$ standard co-ordinate vector of $\bC^{2g+1}$. Using the identification $\tau(\cdot) = \tr_{A_f / \Q}(\cdot / f'(x))$ \cite[Chapter III, \S6]{serre-localfields} of linear forms $A_f \to \Q$, we can compute
    \[ ( e_i, e_i )_{W_D} = \tr_{A_f}(e_i^2 / f'(x) U(x)) = 1 / f'(\omega_i) U(\omega_i). \]
    Therefore $\langle e_i, e_i \rangle_{H_{T_D}} = | f'(\omega_i) U(\omega_i)|^{-1}$, and 
    \[ \langle U, U \rangle_{H_{T_D}} = \sum_{i=1}^{2g+1} \frac{ | U(\omega_i) |^2 }{| f'(\omega_i) U(\omega_i) |} = \sum_{i=1}^{2g+1} \frac{ |U(\omega_i)|}{|f'(\omega_i)|}. \]
    The claimed result follows on using the formula $v_D = N U$. 
\end{proof}

\subsection{The canonical plot of the distinguished orbit}\label{subsec_canonicalplotdistinguishedorbit}

Let $L$ be a lattice (i.e.\ discrete cocompact subgroup) in an inner product space $V$. It is of interest to study not only short vectors in $L$ but `short subgroups', i.e.\ subgroups $L' \leq L$ which have small covolume in the induced inner product space $V' = \R \cdot L' \leq V$. Following Grayson \cite{Grayson-reduction}, we associate to $L$ its \define{canonical plot}, as follows: plot, for each $i = 0, \dots, \rank L$ and for each submodule $L' \leq L$ of rank $i$, the point $(i, \log \covol L')$ in $\R^2$. The canonical plot of $L$ is the polygon starting at $(0, 0)$ and ending at $(\rank L, \log \covol L)$ which is the lower convex hull of these points. It provides a measure of the skewness of $L$: in particular, its slopes are closely related to the Minkowski successive minima of $L$ (see e.g. \cite{Borek-slopes}). Moreover, associated to the canonical plot is the \define{canonical filration} of $L$: let $(0, 0), (i_1, \log \covol(L_1)), \dots, (i_m, \log \covol(L_m)) = (\rank L, \log \covol(L))$ be the vertices of the canonical plot. Then we have 
\[ 0 \leq L_1 \leq L_2 \leq \dots \leq L_m \]
and for each $j = 1, \dots, m$, $L_j$ is the \emph{unique} subgroup of $L$ of rank $i_j$ and covolume $\covol(L_j)$. 

The shape of the canonical plot of $L$ has significant consequences for the presence of other `short subgroups' of $L$ in view of the following lemma:
\begin{lemma}\label{lem_parallelogram_law}
    Let $L_1, L_2 \leq L$ be subgroups. Then we have
    \[ \log \covol(L_1 \cap L_2) + \log \covol(L_1 + L_2) \leq \log \covol(L_1) + \log \covol(L_2). \]
\end{lemma}
\begin{proof}
    See \cite[Theorem 1.12]{Grayson-reduction}.
\end{proof}
Now suppose that $f(x) = x^{2g+1} + c_2 x^{2g-1} + \dots + c_{2g+1} \in \Z[x]$ is a polynomial of nonzero discriminant, and that $(U, V, R)$ is a Mumford triple for $f$ over $\Q$ associated to a divisor $D$. Then we have constructed a $\Z$-lattice $M_D \leq W_D$ satisfying $T_D M_D \leq M_D$. By Proposition \ref{prop_existence_of_reduction_covariant}, $W_D \otimes_\Q \R$ comes equipped with a canonical inner product, compatible with the structure of quadratic space. This compatibility implies that the canonical plot of the lattice $M_D \leq W_D \otimes_\Q \R$ has an additional symmetry: indeed, it is symmetric about the line $x = (2g+1)/2$, because for any saturated subgroup $L' \leq M_D$, we have $\covol L' = \covol (L')^\perp$ (cf. the discussion in \cite[\S 7]{Grayson-reduction}).

The following lemma gives an upper bound on the canonical plot associated to the divisor $D = 0$ (hence to the $\Z$-lattice $M_1 \leq W_1 = A_f$). 
\begin{lemma}
    For each $m = 1, \dots, 2g+1$, let $L_m = \langle 1, x, \dots, x^{m-1} \rangle_\Z$. Then we have
    \[ \log \covol L_m = -\frac{1}{4} \log | \Delta | + \frac{1}{2} \log  \sum_{\substack{J \subset \{ 1, \dots, 2g+1\} \\ |J| = m }} \prod_{\substack{i \in J, j \in J \\ i < j}} | \omega_i - \omega_j| \prod_{\substack{i \not\in J, j \not\in J \\ i < j}} | \omega_i - \omega_j|. \]
\end{lemma}
    It follows from the formula that $\covol L_m = \covol L_{2g+1-m}$. This is a consequence of the fact that $L_m^\perp = L_{2g+1-m}$. 
\begin{proof}
    We compute in a similar way to the proof of Theorem \ref{thm_good_global_lattice}. Let $e_1, \dots, e_{2g+1} \in W_1 \otimes_\Q \bC \cong \bC^{2g+1}$ be the standard co-ordinate vectors. They are pairwise orthogonal, and $\langle e_i, e_i \rangle_{H_{T_1}} = 1 / | f'(\omega_i) |$. The covolume of $L_m$ is equal to the length of the vector $1 \wedge x \wedge \dots \wedge x^{m-1} \in \wedge^m W_D$. Since $x^i = \sum_j \omega_j^i e_j$, the Vandermonde identity shows that 
    \begin{multline*} (\covol L_m)^2 = \sum_J \prod_{\substack{i, j \in J \\ i < j}} | \omega_i - \omega_j|^2 \prod_{j \in J} | f'(\omega_j) |^{-1} \\ = \sum_J \prod_{\substack{i, j \in J \\ i \neq j}} | \omega_i - \omega_j| \prod_{\substack{i \in J, 1 \leq j \leq 2g+1\\i \neq j}} | \omega_i - \omega_j|^{-1} 
    = \sum_J \prod_{\substack{i \in J \\ j \not\in J}} | \omega_i - \omega_j |^{-1}, \end{multline*} 
    where the sum is over all subsets $J \subset \{ 1, \dots, 2g+1 \}$ of size $m$. Pulling out a factor of  $| \Delta |^{-1/2}$, we get
    \[ (\covol L_m)^2 = | \Delta |^{-1/2} \sum_J \prod_{\substack{i \in J, j \in J \\ i < j}} | \omega_i - \omega_j| \prod_{\substack{i \not\in J, j \not\in J \\ i < j}} | \omega_i - \omega_j|. \]
\end{proof}

The next proposition gives a lower bound for the length of a vector in $M_1$ not lying the `short' subgroup $L_g$ and will be useful in giving a lower bound for heights of points in $2J_f(\Q)\setminus \{0\}$.

\begin{proposition}\label{prop_lowerboundvectordistinguishedorbit}
Let $f(x) = x^{2g+1} + c_2 x^{2g-1} + \dots + c_{2g+1} \in \Z[x]$ be a polynomial of nonzero discriminant. Let $\omega_1, \dots, \omega_{2g+1} \in \bC$ be the complex roots of $f$. Suppose give constants $\delta, X > 0$ such that $\height(f) \leq X$ and for each $1 \leq i < j \leq 2g+1$, we have $|\omega_i - \omega_j| \geq X^{1-\delta}$. Then:
    \begin{enumerate}
         \item If $v = \sum_{i=0}^m a_i x^i \in M_1$ is any vector with $g < m < 2g+1$ and $a_m \neq 0$, then we have
            \[ \frac{1}{2} \log \langle v, v \rangle_{H_{T_1}} \geq (1 -  \delta g(2g+1)) \log X - \frac{(2g+1)(4g+1)}{2} \log 2 + \log |a_m|. \]
        \item Suppose that $\delta < 1 / (2g(2g+1))$ and that $\log X > (2g+1)(4g+1) \log 2 / (1 - \delta 2 g (2g+1))$. Then the canonical filtration of $M_1$ equals $0 \leq L_1 \leq L_2 \leq \dots \leq L_{2g} \leq L_{2g+1} = M_1$, where $L_i = \langle 1, x, \dots, x^{i-1} \rangle_\Z$. 
    \end{enumerate}
\end{proposition}
    \begin{proof}
        Let $N = 2g(2g+1)$. We first bound $\log \covol L_i$ for each $i = 1, \dots, 2g$. By Lemma \ref{lem_height_comparison}, we have $\max_i | \omega_i | \leq 2 X$, hence $X^{N(1-\delta)} \leq \Delta \leq (4X)^N$, and
        \begin{multline*} \log \covol L_i \leq -\frac{1}{4} N(1-\delta) \log X + \frac{1}{2} \log (2^{(2g+1)^2} X^{g(2g+1) +i(i-2g-1)}) \\
        = \frac{(2g+1)^2}{2} \log 2 + \frac{1}{4} \left( \delta N + 2i(i-2g-1 ) \right) \log X.
        \end{multline*} 
        In the other direction, we have
        \begin{multline*} \log \covol L_i \geq -\frac{1}{4} N \log 4X + \frac{1}{2} \log X^{(g(2g+1) + i(i-2g-1))(1-\delta)} \\ 
        =  -g(2g+1) \log 2 + \frac{1}{4}( -\delta N + 2i(i-2g-1)) \log X.
        \end{multline*} 
        In particular, we find      \begin{equation}\label{eqn_difference_in_covolumes} \log \covol L_{i+1} - \log \covol L_i \geq (i - g - \delta N / 2) \log X - \frac{(2g+1)(4g+1)}{2} \log 2. 
        \end{equation}
        We can now prove the first part of the proposition. Let $v = \sum_{i=0}^m a_i x^i \in M_1$ be as in the statement. We apply Lemma \ref{lem_parallelogram_law} to the subgroups $L_m, \langle v \rangle_\Z \leq M_1$. We find that $L_m \cap \langle v \rangle_\Z = 0$, while $L_m + v$ is a subgroup of $L_{m+1}$ of index $|a_m|$, and therefore
        \[ \frac{1}{2} \log \langle v, v \rangle_{H_{T_1}} + \log \covol L_m \geq \log |a_m| + \log \covol L_{m+1}. \]
        Re-arranging and using (\ref{eqn_difference_in_covolumes}) gives a formula
        \[ \frac{1}{2} \log \langle v, v \rangle_{H_{T_1}} \geq \log |a_m| + (m-g - \delta N / 2) \log X - \frac{(2g+1)(4g+1)}{2} \log 2.  \]
        Since $m - g$ is a positive integer, this implies the desired bound. 

        To prove the second part, we apply \cite[Corollary 1.30]{Grayson-reduction}, which says that $0 \leq L_1 \leq \dots \leq L_{2g} \leq L_{2g+1} =  M_1$ will be the canonical filtration as soon as the sequence $\log \covol L_{i+1} - \log \covol L_i$ is strictly increasing. We have   \begin{equation}\label{eqn_difference_in_covolumes_lower_bound}  \log \covol L_{i+1} - \log \covol L_i \leq  ( i - g + \delta N / 2) \log X + \frac{(2g+1)(4g+1)}{2} \log 2 . \end{equation}
        Subtracting (\ref{eqn_difference_in_covolumes_lower_bound}) from (\ref{eqn_difference_in_covolumes}), we find
        \[ (\log \covol L_{i+1} - \log \covol L_i) - (\log \covol L_i - \log \covol L_{i-1}) \geq (1- \delta N ) \log X - (2g+1)(4g+1) \log 2. \]
        The right-hand side does not depend on $i$; as soon as it is positive, the cited result applies, leading to the desired result.  
    \end{proof}

\subsection{Further study of integral distinguished orbits}\label{subsec_Qinvariant}
    In this section we record some properties of the $|Q|$-invariant of integral vectors in the distinguished orbit, which was defined and studied in \cite{Bhargava-squarefree}. The conventions of that paper are slightly different: they consider the action of the group $G$ on the space $V'$ of symmetric matrices (by $\gamma \cdot X = \gamma X ({}^t \gamma)$). The embedding $V \to V'$ given by $T \mapsto T J$ intertwines the two $G$-actions. 
    \begin{definition}
        Let $T \in V(\Z)$ be a self-adjoint linear operator of nonzero discriminant, and let $L \leq W(\Q)$ be a maximal isotropic subspace of dimension $g$ such that $T L \leq L^\perp$ (so that $T$ is in the distinguished orbit). Choose $h \in G(\Z)$ such that $h L = \langle e_g, \dots, e_1 \rangle_\Q$, so that $h T h^{-1}$ has the form
        \[  h T h^{-1} = \left( \begin{array}{ccc} A_{11} & 0 \\ A_{21} & A_{22} \end{array}\right) \]
        for matrices $A_{11} \in \Mat_{g \times g+1}(\Z)$, $A_{21} \in \Mat_{g+1 \times g+1}(\Z)$, $A_{21} \in \Mat_{g+1 \times g}(\Z)$. Then $h T h^{-1} J$ lies in the subspace $W_0$ of $V'$ defined on \cite[p. 1045]{Bhargava-squarefree}, and we define $|Q|(T, L)$ to be the absolute value of the number $Q(h T h^{-1} J)$ defined by \cite[Eqn. (7)]{Bhargava-squarefree}.
    \end{definition}
    It follows from \cite[Eqn. (8)]{Bhargava-squarefree} that $|Q|(T, L)$ is a positive integer which is independent of the choice of matrix $h$. When the characteristic polynomial of $T$ is irreducible, there is (cf. \cite[Proposition 2.1]{Bhargava-squarefree}) a unique choice of subspace $L$, leading to the following definition.
    \begin{definition}
        Let $T \in V(\Z)$ be a self-adjoint linear operator with irreducible characteristic polynomial, and suppose that $T$ is in the distinguished orbit. Then we define $|Q|(T) = |Q(T, L)|$, where $L \leq W(\Q)$ is the unique maximal isotropic subspace such that $T L \leq L^\perp$.
    \end{definition}
    A key result in \cite{Bhargava-squarefree} is that irreducible polynomials $f(x)$ admitting a vector in the distinguished orbit of large $|Q|$-invariant are rare, in some sense. Proposition \ref{prop_construction_of_vectors_of_large_Q_invariant} below will show that such vectors can be constructed from lattices in $W_1$ which are far from the standard one $M_1 = \Z[x] / (f(x)) \leq W_1 = \Q[x] / (f(x))$.
    \begin{lemma}\label{lem_change_in_Q_invariant}
        Let $T \in V(\Z)$ be a self-adjoint linear operator of nonzero discriminant, and let $L \leq W(\Q)$ be a maximal isotropic subspace such that $T L \leq L^\perp$. Suppose given $h \in G(\Q)$ such that $h L \leq L$ and $h T h^{-1} \in V(\Z)$. Let $h_L : L \to L$ denote the restriction of $h$ to $L$. Then $|Q|(h T h^{-1}, L) = |\det h_L|^{-1} \times |Q|(T, L)$.
    \end{lemma}
    \begin{proof}
    After first acting by an element of $G(\Z)$, we can assume that $L = \langle e_g, \dots, e_1 \rangle_\Q$ and that $h$ has the block form
    \[ h = \left( \begin{array}{ccc} h_1 & 0 & 0 \\ \ast & 1 & 0 \\ \ast & \ast & h_2 \end{array}\right), \]
    where $h_1, h_2$ are $g \times g$ matrices with $\det h_1 \times \det h_2 = 1$. Then \cite[Eqn. (8)]{Bhargava-squarefree} shows that $|Q|(h T h^{-1}, L) = | \det h_1 | \times |Q|(T, L) = | \det h_2 |^{-1} \times |Q|(T, L)$, as required. 
    \end{proof}
\begin{proposition}\label{prop_construction_of_vectors_of_large_Q_invariant}
    Let $f(x) = x^{2g+1} + c_2 x^{2g-1} + \dots c_{2g+1} \in \Z[x]$ be a polynomial of nonzero discriminant, and let $M' \leq W_1$ be a $T_1$-invariant $\Z$-lattice. Suppose we can find an integer $N > 1$ and a vector $v \in M'$ which is saturated in $\frac{1}{N} M_1$, in the sense that $v \in \frac{1}{N} M_1$ but $v \not\in \frac{p}{N} M_1$ for any prime $p | N$. Then we can find $T_v \in W(\Z)$ of characteristic polynomial $f(x)$ and a maximal isotropic subspace $L \leq W(\Q)$ such that $T_v L \leq L^\perp$ and $N$ divides $|Q|(T_v, L)$. 
\end{proposition}
\begin{proof}
    We consider the flag $F_i = \langle 1, x, \dots, x^{i-1} \rangle_\Q \leq W_1$. By Proposition \ref{prop_uniqueness_of_cusp}, we can find a $\Z$-basis $b_{2g+1}, \dots, b_1$ of $W_1$ such that $\langle b_i, b_j \rangle = \delta_{i, 2g+2-j}$ and $F_i = \langle b_i, \dots, b_1 \rangle_\Q$. We can write uniquely $b_i = \frac{1}{n_i} f_i(x)$ for some positive rational number $n_i$ and monic polynomial $f_i(x) \in \Q[x]$ of degree $i-1$. We claim that in fact we have:
    \begin{itemize}
        \item $n_i n_{2g+2-i} = 1$ for each $i = 1, \dots, 2g+1$. In particular, $n_{g+1} = 1$.
        \item Let $N_i = F_i \cap M_1$. Then $n_i N_i \leq M_1$ ($i = 1, \dots, 2g+1$). In particular, $f_i(x) \in \Z[x]$.
        \item We have $n_{g+2}, n_{g+3}, \dots, n_{2g+1} \in \Z$ and $n_{g+2} | n_{g+3} | \dots | n_{2g+1}$.
    \end{itemize}
    The first point is immediate since $n_i n_{2g+2-i} = ( b_i, b_{2g+2-i} )_{W_1} = 1$. For the second, we use induction on $i$, the case $i = 1$ being immediate (since $b_1 = n_1^{-1}$ is constant). In general, we observe $M'$ is $x$-stable, so $x b_i \in N_{i+1}$ and $x b_i \equiv \frac{n_{i+1}}{n_i} b_{i+1} \text{ mod } F_i$, hence $n_{i+1} \in n_i \Z$. We have $\frac{n_{i+1}}{n_i} b_{i+1} - x b_i \in F_i \cap M' = N_i$, hence $n_{i+1} b_{i+1} - n_i x b_i \in M_1$, hence $f_{i+1}(x) \in M_1$ and $n_{i+1} N_{i+1} \leq M_1$. This proves the second point. Since $n_{g+1} = 1$ and $n_{i+1} \in n_i \Z$ for each $i = 1, \dots, 2g$, the third point now also follows. 

    Now let $v \in M'$ be the vector in the statement of the proposition. We claim that $N | n_{2g+1}$. By assumption, we can write $v = \sum_{i=1}^{2g+1} a_i b_i$ for some integers $a_i$. Let $p$ be a prime such that $v_p(N) = l > 0$. It is enough to show that $v_p(n_{2g+1}) \geq l$. We have $v = \sum_{i = 1}^{2g+1} \frac{a_i}{n_{i}} f_i(x)$ with $f_i(x) \in M_1$. Suppose that $v_p(n_{2g+1}) < l$. Since $v_p(n_{2g+1}) \geq v_p(n_{2g}) \geq \dots \geq v_p(n_{g+1}) = 0$, this implies that $v \in p^{1-l} M_1$, a contradiction.  
    
    With this claim in hand, we can now complete the proof. Let $\beta_{2g+1}, \dots, \beta_{1}$ be a $\Z$-basis of $W_1$ such that $( \beta_i, \beta_j)_{W_1} = \delta_{i, 2g+2-j}$ and $F_i = \langle \beta_1, \dots, \beta_{i} \rangle_\Q$. Let $\varphi : W_1 \to W(\Q)$ be the isomorphism of quadratic spaces which sends $\beta_i$ to $e_i$ ($i = 1, \dots, g$), $\beta_{g+1}$ to $e_0$, and $\beta_{2g+2-i}$ to $e_{-i}$ ($i = g, \dots, 1$). Then $\varphi(M_1) = W(\Z)$, $\varphi(F_g) = \langle e_g, \dots, e_{1} \rangle_\Q = L$, say, and if we set $T = \varphi \circ T_1 \circ \varphi^{-1}$, then $T \in V(\Z)$ and $|Q|(T, \varphi(F_g)) = 1$ (as follows by direct computation).

    Let  $\psi : W_1 \to W(\Q)$ be the isomorphism of quadratic spaces that sends $b_i$ to $e_i$ ($i = 1, \dots, g$), $b_{g+1}$ to $e_0$, and $b_{2g+2-i}$ to $e_{-i}$ ($i = g, \dots, 1$), and let $h = \psi \circ \varphi^{-1} \in G(\Q)$. Then $T_v = \psi \circ T_1 \circ \psi^{-1} \in V_f(\Z)$ and $T_v L \leq L^\perp$. The matrix $h$ is lower-triangular, with diagonal entries $(n_{2g+1}, n_{2g}, \dots, n_1)$, and satisfies $T_v = h \cdot T$. 

    Applying Lemma \ref{lem_change_in_Q_invariant}, we finally conclude that
    \[ |Q|(T_v, L) = (n_1 \dots n_g)^{-1} |Q|(T, L) = n_{g+2} \dots n_{2g+1}. \]
    Since $n_{g+2}, \dots, n_{2g+1}$ are integers, and since $N | n_{2g+1}$, we find that $N$ divides $|Q|(T_v, L)$, as required. 
\end{proof}

\section{A density 1 family}\label{sec_density_1_family}

If $\delta, X > 0$, then we define the following sets of polynomials $f(x) = x^{2g+1} + c_2 x^{2g-1} + \dots + c_{2g+1} \in \Z[x]$:
\begin{itemize}
    \item $\mathcal{F}(X)$ denotes the set of polynomials $f(x)$ such that $\Delta(f) \neq 0$ and $\height(f) < X$.
    \item $\mathcal{M}_{\delta, 1}(X)$ denotes the set of $f(x) \in \mathcal{F}(X)$ such that  for all $i\neq  j$ we have $|\omega_i - \omega_j| > X^{1-\delta}$ (where $\omega_1, \dots, \omega_{2g+1}$ are the complex roots of $f$).
    \item $\mathcal{M}_{\delta, 2}(X)$ denotes the set of $f(x) \in \mathcal{F}(X)$ such that $f(x)$ is irreducible and for each element $T \in V(\Z)$ lying in the distinguished orbit of characteristic polynomial $f$, we have $|Q|(T) < X^\delta$, in the notation of \S\ref{subsec_Qinvariant}.
    \item $\mathcal{F}_\delta(X) = \mathcal{M}_{\delta, 1}(X) \cap \mathcal{M}_{\delta, 2}(X)$.
\end{itemize}

Recall that $\#\mathcal{F}(X) = 2^{2g} X^{g(2g+3)} + O(X^{g(2g+3)-1})$ by Lemma \ref{lemma_F(X)_expectedsize}.

\begin{proposition}\label{prop_F_delta_has_density_1}
    For any $\delta > 0$, we have $\lim_{X \to \infty} \# \mathcal{F}_\delta(X) / \# \mathcal{F}(X) = 1$.
\end{proposition}
\begin{proof}
    We will show that $\# (\mathcal{F}(X) - \mathcal{M}_{i, \delta}(X)) = o(X^{g(2g+3)})$ for $i = 1, 2$.
    For $i = 1$, we observe that $
    \mathcal{M}^c_{1, \delta}(X) = \mathcal{F}(X) - \mathcal{M}_{1, \delta}(X)$ is the set of polynomials such that $\height(f) < X$ and $| \omega_i - \omega_j | \leq X^{1-\delta}$ for some $i \neq j$. In particular, if $f(x)$ lies in this set then $| \Delta(f) | \ll X^{2g(2g+1) - 2 \delta}$. Let $U_\delta(X) \leq \R^{2g}$ denote the set of vectors $(c_2, \dots, c_{2g+1})$ such that $|c_i| \leq X^{1/i}$ and the polynomial $f(x)$ has two roots $\omega_i, \omega_j$ such that $|\omega_i - \omega_j| \leq X^{1-\delta}$ By Proposition \ref{prop_countlatticepointsbarroero}, it suffices to show that $\vol(U_\delta(X)) = o(X^{g(2g+3)})$ as $X \to \infty$. For each $0 \leq m \leq g$, let $V_{\delta, m}(X) \subset \bC^{2m} \times \R^{2g+1-2m} \subset \bC^{2g+1}$ denote the set of tuples $(\omega_i)$ such that $\omega_i = \overline{\omega}_{m+i}$ ($i = 1, \dots, m$), $\sum_{i=1}^{2g+1} \omega_i = 0$, and $| \omega_i | \leq 2 X$ for each $i$, and $|\omega_i - \omega_j| \leq X^{1-\delta}$ for some $i \neq j$. Then $U_\delta(X)$ is covered by the images of the sets $V_{\delta, m}(X)$, and by the Jacobian change of variable formula, we have
    \[ \vol(U_\delta(X)) \leq \sum_m \int_{\omega \in V_{\delta, m}} \prod_{1 \leq i < j \leq 2g+1} |\omega_i - \omega_j|\,d \omega \ll \sum_m X^{g(2g+1) - \delta} \vol(V_{\delta, m}(X)). \]
    It therefore suffices to show that we have $\vol(V_{\delta, m}(X)) \ll X^{2g}$ for each $m$, and this is clear. 

    For $i = 2$, we already know that the number of elements of $\mathcal{F}(X)$ which are reducible is $o(X^{g(2g+3)})$ (cf. \cite[Lemma 4.3]{Bhargava-squarefree}). It therefore suffices to prove that the number of elements $f(x) \in \mathcal{F}(X)$ which are irreducible and admit a distinguished integral vector $T \in V(\Z)$ of characteristic polynomial $f$ and $|Q|(T) > X^\delta$ is $o(X^{g(2g+3)})$. This number is bounded above by the quantity $N(\mathcal{L}, M, X)$ introduced on \cite[p. 1050]{Bhargava-squarefree} (roughly, the number of orbits with $Q$-invariant $>M$ and height $<X$), with $M$ taken to be $X^\delta$. Now \cite[Propositions 2.5, 2.6, 2.7]{Bhargava-squarefree} together show that $N(\mathcal{L}, M, X) = o(X^{g(2g+1)})$, even with a power-saving. (In fact, one needs to restrict the statements here to the trace-0 subspace of the relevant representation, but as observed in the proof of \cite[Theorem 5.4]{Bhargava-squarefree}, the arguments go through essentially unchanged in this case.)
\end{proof}
The following statement was referred to in the introduction, and is not used elsewhere in this paper. 
\begin{proposition}\label{prop_statistical_Lang}
    Suppose that $g = 1$, and let $\epsilon > 0$. Then
    \[ \lim_{X \to \infty} \frac{\# \{ f \in \mathcal{F}(X) \mid \exists P \in C_f(\Q), \hat{h}(P) \leq (\frac{1}{2} - \epsilon) \log \height(f)\} }{ \# \mathcal{F}(X) } = 0, \]
    where $\hat{h} : C_f(\Q) \to \R$ denotes the usual canonical height of the elliptic curve $C_f$ with origin $P_\infty$.
\end{proposition}
\begin{proof}
    In this proof, the symbol $d_i$ ($i = 0, 1, \dots$) will denote an absolute constant. Let $\mathcal{M}_\delta(X)$ denote the set of polynomials $f(x) = x^3 + c_2 x + c_3$ denoted by $\mathcal{F}_\delta(X)$ in \cite[\S 2]{LeBoudec-statistical} (note that the normalisations there are slightly different -- in particular, we have $\height(f) \ll X^{1/6}$ if $f \in \mathcal{M}_\delta(X)$). We follow the proof of \cite[Lemma 2]{LeBoudec-statistical}. Let $f(x) \in \mathcal{M}_\delta(X)$. We use the decomposition $\hat{h} = \sum_p \lambda_p + \lambda_\infty$ as a sum of local heights. Summing \cite[(2.15), (2.16)]{LeBoudec-statistical}, we find that if $P \in C_f(\Q) - \{ P_\infty \}$ then
    \[ \sum_p \lambda_p(P) \geq \frac{1}{12} \log | \Delta | - \frac{3}{8} \delta \log X + d_1. \]
     Using \cite[Proposition 5.4]{Silverman-differenceweilcanonical} and the bound $|j(C_f)| \leq 1728 \cdot 16 X^\delta$, we find
    \[ \lambda_\infty(P) \geq - \frac{1}{8} \delta \log X + d_2. \]
     Summing up gives 
    \[ \hat{h}(P)\geq \frac{1}{12} \log | \Delta | - \frac{1}{2} \delta \log X  + d_3 \geq \frac{1}{12} (1 - 7 \delta) \log X + d_3 \]
    Since $f \in \mathcal{M}_\delta(X)$, it follows from the definition that $\height(f) \leq \frac{1}{6} \log X + d_4$. Since (\cite[Lemma 1]{LeBoudec-statistical}) $\mathcal{M}_\delta(X)$ is a density 1 family, the claimed result now follows easily on choosing $\delta$ sufficiently small compared to $\epsilon$. 
\end{proof}
    
\section{Lower bounds for height functions}\label{sec_lower_bounds_for_height_functions}

In this section, we will prove our main theorems.

\subsection{Height functions on $J_f$ and statement of main theorem}

Given a point $P \in \P^n(\Q)$, let $h(P)\in \R_{\geq 0}$ denote the logarithmic Weil height of $P$, defined by the formula
\[ h( [ x_0 : x_1 : \dots : x_n ]) = \sum_v \log \max( |x_0|_v, |x_1|_v, \dots, |x_n|_v), \]
where $v$ runs over the set of places of $\Q$, and the absolute values $|\cdot|_v$ are normalised in such a way that they satisfy the product formula. If $\alpha\in \Q$, we set $h(\alpha) \coloneqq h([\alpha:1])$.

Let $f = x^{2g+1} +c_2x^{2g-1} + \cdots + c_{2g+1}\in \Q[x]$ be a separable polynomial and $J_f$ the Jacobian of the hyperelliptic curve $C_f$. Pulling back the theta divisor on $\Pic^{g-1}_{C_f}$ along $(g-1)P_{\infty}$ defines a line bundle $\mathcal{L}_f$ corresponding to the principal polarization on $J_f$. Then (\cite[Theorem 4.8.1]{BirkenhakeLangeAV}) $\mathcal{L}_f^{\otimes 2}$ defines an embedding of the Kummer variety $J_f / \{ \pm 1 \}$ in $\P(H^0(J_f, \mathcal{L}_f^{\otimes 2}))$, and one way to define a naive height on $J_f$ is to choose a basis for the space of sections, leading to a morphism $\varphi_f : J_f \to \P^{2^g-1}$, and pull back the height from projective space. Let us write $h_{\mathcal{L}_f}(P) = \frac{1}{2} h( \varphi_f(P) ) : J_f(\Q) \to \R_{\geq 0}$ for the height defined this way, which depends on the choice of co-ordinates. 

In this paper, we use another height, defined in terms of the Mumford representation of a hyperelliptic curve, that makes sense in families, and does not require any additional choices. This height has already been studied in \cite{holmes-arakelovapproachnaiveheighthyperelliptic}.

Recall from \S\ref{subsec_rationalorbits} that every $[D]\in J_f(\Q)$ corresponds to a Mumford triple $(U,V,R)$ where $U = x^m + u_1 x^{m-1} +\dots + u_m \in \Q[x]$ is monic of degree $m\leq g$.
We define
\begin{align}
    h^{\dagger}([D]) \coloneqq h([1:u_1:\dots : u_m]).
\end{align}
If $g = 1$ and $D = [\alpha : \beta : 1 ] - P_\infty$ for a rational point $(\alpha, \beta) \in C_f^0(\Q)$, then $h^\dagger([D]) = h(\alpha)$ is the usual naive height associated to an elliptic curve in short Weierstrass form.

\begin{proposition}\label{prop_lower_bound_for_usual_naive_height}
    We have $h_{\mathcal{L}_f} \geq \frac{1}{2} h^\dagger + O_f(1)$ on $J_f(\Q)$.
\end{proposition}
Although not necessary for the statement of our main results, Proposition \ref{prop_lower_bound_for_usual_naive_height} shows that a lower bound for $h^\dagger$ is meaningful in terms of the usual naive height on $J_f(\Q)$.
\begin{proof}
    This follows from \cite[Theorem 37]{holmes-arakelovapproachnaiveheighthyperelliptic} (which gives the statement for the canonical height $\hat{h}_{\mathcal{L}_f}$ -- but $|\hat{h}_{\mathcal{L}_f} - h_{\mathcal{L}_f}|$ is bounded on $J_f(\Q)$).
\end{proof}
In some cases we can be more precise. For example, in the case $g = 2$, Cassels and Flynn \cite[Ch. 3]{CasselsFlynnProlegomena} write down equations for the Kummer surface in $\P^3$ with respect to a fixed choice of co-ordinates, depending only on the coefficients of $f(x)$. One then obtains a height $h_{\mathcal{L}_f}$ on $J_f(\Q)$ not depending on any further choices. The equations for the morphism $J_f \to \P^3$ given in \emph{loc. cit.} make it clear that we have $h_{\mathcal{L}_f}(P) \geq \frac{1}{2} h^\dagger(P)$ for all $P \in J_f(\Q)$ (and not just up to a constant depending on $f$). 

We can now state our main theorem.
\begin{theorem}\label{thm_main_theorem_on_heights}
    Let $\epsilon > 0$. Then we have
    \[ \lim_{X \to \infty} \frac{ \# \{ f \in \mathcal{F}(X) \mid \forall P \in J_f(\Q)\setminus\{0\}, h^\dagger(P) \geq (g - \epsilon) \log \height(f) \}}{\# \mathcal{F}(X) } = 1. \]
\end{theorem}
The proof of this theorem will occupy the remainder of \S \ref{sec_lower_bounds_for_height_functions}.
The proof of Theorem \ref{thm_intro_degree1} follows from a slight modification of the argument and is omitted.

\subsection{Proof of Theorem \ref{thm_main_theorem_on_heights}}

We begin with some preparatory lemmas, which bring together all of the ideas developed so far. 
\begin{lemma}\label{lemma_upper_bound_for_norm_of_U}
    Let $f \in \mathcal{F}_\delta(X)$, and let $U(x)= x^m + u_1 x^{m-1} + \dots + u_m \in \Q[x]$ be a polynomial of degree $1 \leq m \leq g$ and with no roots in common with $f$. Then we have
    \[ \frac{1}{2} \log \sum_{i=1}^{2g+1} \frac{|U(\omega_i)|}{|f'(\omega_i)|} \leq \frac{1}{2}(m - 2g + \delta g(2g+1)) \log X + \frac{1}{2} \log \max(1, |u_1|, \dots, |u_m|) + c_{g, m},  \]
    where $c_{g, m} = \frac{1}{2}\log[(m+1)(2g+1) 2^{m+2g(2g-1)}]$. In particular, if $M \geq 1$ is the least natural number such that $M U \in \Z[x]$, then
    \[ \frac{1}{2} \log M + \frac{1}{2} \log \sum_{i=1}^{2g+1} \frac{|U(\omega_i)|}{|f'(\omega_i)|} \leq \frac{1}{2} h([1 : u_1 : \dots : u_m]) + \frac{1}{2}(m - 2g + \delta g(2g+1)) \log X + c_g. \]
\end{lemma}
\begin{proof}
    This is an elementary computation.
    Let $\omega_1, \dots, \omega_{2g+1}$ be the complex roots of $f$. Since $f \in \mathcal{F}_\delta(X)$, we have $X^{1-\delta} \leq |\omega_i - \omega_j| \leq 4X$ for all $1 \leq i < j \leq 2g+1$, hence $|\Delta| \geq X^{2g(2g+1)(1-\delta)}$. We have
    \begin{align*}  
    \sum_{i=1}^{2g+1} \frac{|U(\omega_i)|}{|f'(\omega_i)|} &= |\Delta|^{-1/2} \sum_{i=1}^{2g+1} |U(\omega_i)| \prod_{\substack{j < k \\ j, k \neq i}} | \omega_j - \omega_j|\\
    &\leq X^{-g(2g+1)(1-\delta)} (4X)^{2g(2g-1)/2}\sum_{i=1}^{2g+1} |U(\omega_i)| \\
    &= 2^{2g(2g-1)} X^{-2g + \delta g(2g+1)} \sum_{i=1}^{2g+1} |U(\omega_i)|.
    \end{align*} 
    We also have
    \[ | U(\omega_i)| \leq \max(1, |u_1|, \dots, |u_m|) (1 + |\omega_i| + \dots + |\omega_i|^{m}) \leq (m+1) (2X)^m \max(1, |u_1|, \dots, |u_m|)   \]
    (noting that the existence of $f$ implies that $X \geq 1$), hence
    \[ \sum_{i=1}^{2g+1} | U(\omega_i)| \leq (2g+1)(m+1) 4^m X^m \max(1, |u_1|, \dots, |u_m|). \]
    Putting these estimates together gives
    \begin{multline*}  \frac{1}{2} \log \sum_{i=1}^{2g+1} \frac{|U(\omega_i)|}{|f'(\omega_i)|} \leq \\ (m/2 - g + \delta g(2g+1)/2) \log X + \frac{1}{2} \log \max(1, |u_1|, \dots, |u_m|) + \frac{1}{2}\log[(m+1)(2g+1) 2^{m+2g(2g-1)}], 
    \end{multline*}
    as required. 
\end{proof}

The next lemma quantifies the idea that if $P\in J_f(\Q)$ has small height, then there exists a $T\in V_f(\Z)$ and a vector $v\in W(\Z)$ of small height.

\begin{lemma}\label{lem_height_of_non_divisible_point}
    Let $f \in \mathcal{F}_\delta(X)$, and let $P \in J_f(\Q)$ be nontrivial, of Mumford degree $m$. Then there exists a vector $T \in V_f(\Z)$ with the following properties:
    \begin{enumerate}
        \item $T$ lies in the rational orbit associated to $P$ (in the sense of Corollary \ref{corollary_rationalorbits_exist}). In particular, if $P \not\in 2 J_f(\Q)$ then $T$ is irreducible.
        \item There exists a nonzero vector $v \in W(\Z)$ such that 
        \[ \frac{1}{2} h^\dagger(P) \geq \frac{1}{2} \log \langle v, v \rangle_{H_T} + \frac{1}{2}(2g - m - \delta g(2g+1)) \log X - c_{g, m}. \]
    \end{enumerate}
\end{lemma}
\begin{proof}
    Let $(U, V, R)$ be the Mumford triple associated to $P$. Then $0 < \deg U \leq g$; since $f \in \mathcal{F}_\delta(X)$, $f(x)$ is irreducible and in particular $U$, $f$ have no roots in common. Let $M \geq 1$ be the least integer such that $M U \in \Z[x]$. By Theorem \ref{thm_good_global_lattice} and Proposition \ref{prop_classification_integralorbits}, we can find $T \in V_f(\Z)$ and a primitive vector $v \in W(\Z)$ such that 
    \[ \frac{1}{2} \log \langle v, v \rangle_{H_T}  = \frac{1}{2} \log M + \frac{1}{2} \log \sum_{i=1}^{2g+1} \frac{ |U(\omega_i)|}{|f'(\omega_i)|}.  \]
    The result follows from this equation and Lemma \ref{lemma_upper_bound_for_norm_of_U}. 
\end{proof}

The next lemma studies points in $2J_f(\Q)\setminus \{0\}$ using the results of \S\ref{subsec_canonicalplotdistinguishedorbit} and \S\ref{subsec_Qinvariant}.

\begin{lemma}\label{lem_height_of_divisible_point_non_reduced_case}
    Let $f \in \mathcal{F}_\delta(X)$, and let $P = [D] \in 2 J_f(\Q)$ be nontrivial of Mumford degree $m$. Then either:
    \begin{enumerate}
        \item There exists a reduced divisor $D_1$ such that $D = 2 D_1$ (equality of divisors);
        \item Or we have
        \[ \frac{1}{2} h^\dagger(P) \geq (1 + g - m/2 - \delta(3g(2g+1)/2 + 1)) \log X - c_{g, m}', \]
        where $c_{g, m}' = c_{g, m} + \log 2^{(2g+1)(4g+1)/2}$.
    \end{enumerate}
\end{lemma}
\begin{proof}
    Let $(U, V, R)$ be the Mumford triple associated to $P$. Since $P \in 2 J_f(\Q)$, the class $U \text{ mod }f \in A_f^\times$ is a square; therefore we can write $U \equiv U_1^2 \text{ mod }f$ for a unique polynomial $U_1 \in \Q[x]$ of degree $\leq 2g$. If $\deg U_1 \leq g$, then the congruence $U \equiv U_1^2 \text{ mod }f$ implies that in fact $U = U_1^2$. Since the divisor $D$ is reduced, this implies that each point of $C_f$ must occur in $D$ with even multiplicity and therefore that $D = 2 D_1$ for a uniquely determined reduced divisor $D_1$. 

    Suppose instead that $\deg U_1 > g$. Let $M \geq 1$ be the least integer such that $M U \in \Z[x]$.
    By Theorem \ref{thm_good_global_lattice}, $M = N^2$ is a square and there exists a $T_D$-invariant lattice $M_D \leq W_D$ such that $N U \in M_D$ is primitive. We consider the isomorphism $W_D \to W_1$ of quadratic spaces given by multiplication by $U_1^{-1}$. Let $M'$ denote the image of $M_D$ in $W_1$ under this isomorphism, and $v'$ the image of $N U$ in $M'$; then in fact $v' = N U_1$, and in particular $\deg v' > g$. 

    Since $f \in \mathcal{F}_\delta(X)$, there are no vectors in $V_f(\Z)$ of $|Q|$-invariant $> X^\delta$. By Proposition \ref{prop_construction_of_vectors_of_large_Q_invariant}, there exists an integer $1 \leq \alpha \leq X^\delta$ such that $\alpha v' \in M_1$. We can then apply Proposition \ref{prop_lowerboundvectordistinguishedorbit} to conclude that
    \[ \frac{1}{2} \log \langle \alpha v', \alpha v' \rangle_{H_T} \geq (1 - \delta g(2g+1)) \log X - \log 2^{(2g+1)(4g+1)/2}. \]
    Since $\alpha v' = N U_1$ and $\alpha \leq X^\delta$, we have
    \[ \frac{1}{2} \log \langle \alpha v' \alpha v' \rangle_{H_T} \leq \delta \log X + \frac{1}{2} \log M + \frac{1}{2} \log \sum_{i=1}^{2g+1} \frac{| U(\omega_i) | }{|f'(\omega_i)|}. \]
    We obtain
    \[ \frac{1}{2} \log M + \frac{1}{2} \log \sum_{i=1}^{2g+1} \frac{| U(\omega_i) | }{|f'(\omega_i)|} \geq (1 - \delta(g(2g+1) + 1)) \log X - \log 2^{(2g+1)(4g+1)/2}. \]
    The result now follows on combining this estimate with Lemma \ref{lemma_upper_bound_for_norm_of_U}. 
\end{proof}
\begin{lemma}\label{lem_height_of_divisible_point_reduced_case}
    Let $D$ be a reduced divisor of Mumford degree $m$ such that $2D$ is also reduced. Then $h^\dagger(2D) \geq 2 h^\dagger(D) - \log 2^{6m-2}$.
\end{lemma}
\begin{proof}
    Let $(U, V, R)$ be the Mumford representation of $D$. Then $U^2$ is the first part of the Mumford representation of $2D$, and we need to compare the height of $U$ to the height of $U^2$. The desired inequality follows immediately from \cite[Ch. VIII, Th. 5.9]{Silverman-arithmeticellcurvesbook}. 
\end{proof}
We can now complete the proof of Theorem \ref{thm_main_theorem_on_heights}. 
\begin{proof}[Proof of Theorem \ref{thm_main_theorem_on_heights}]
Let $\epsilon > 0$. We must show that 
\[ \lim_{X \to \infty} \frac{ \# \{ f \in \mathcal{F}(X) \mid \exists P \in J_f(\Q) - \{ 0 \}, \frac{1}{2} h^\dagger(P) \leq (g/2 - \epsilon) \log \height(f) \} }{ \# \mathcal{F}(X) } = 0. \]
By Proposition \ref{prop_F_delta_has_density_1}, it will suffice to show that for some $\delta > 0$, we have
\begin{align}\label{eq_smallpointsdensityzero}
    \lim_{X \to \infty} \frac{ \# \{ f \in \mathcal{F}_\delta(X) \mid \exists P \in J_f(\Q) - \{ 0 \}, \frac{1}{2} h^\dagger(P) \leq (g/2 - \epsilon) \log X \} }{ \# \mathcal{F}(X) } = 0. 
\end{align}
Choose $\delta$ such that $\delta (3g(2g+1)/2 + 1) < \epsilon$ and let $\mathcal{F}(X)^{\text{bad}}$ be the subset of $\mathcal{F}_{\delta}(X)$ appearing in the numerator of \eqref{eq_smallpointsdensityzero}.
Lemma \ref{lem_height_of_divisible_point_non_reduced_case} and Lemma \ref{lem_height_of_divisible_point_reduced_case} together show that there is $X_0 > 0$ such that if $X > X_0$ and $f \in \mathcal{F}(X)^{\text{bad}}$, then there is $P \in J_f(\Q)- 2J_f(\Q)$ such that $\frac{1}{2} h^\dagger(P) < (g/2 - \epsilon) \log X$. 

However, Lemma \ref{lem_height_of_non_divisible_point} shows, after possibly increasing $X_0$, that if $X > X_0$ and $f \in \mathcal{F}_\delta(X)$ has the property that there is a $P \in J_f(\Q) - 2 J_f(\Q)$ with $\frac{1}{2} h^\dagger(P) < (g/2 - \epsilon) \log X$, then there is an endomorphism $T \in V_f(\Z)$, not in the distinguished orbit, and a vector $v \in W(\Z)$ such that $\frac{1}{2} \log \langle v, v \rangle_{H_T} < ( \delta g (2g+1)/2 - \epsilon) \log X$. 

To complete the proof, it therefore suffices to show that for any $\epsilon_1 > 0$, we have
\[ \lim_{X \to \infty} \frac{ \# \{ f \in \mathcal{F}(X) \mid \exists T \in V_f(\Z)^{\text{irr}}, w \in W(\Z) \text{ s.t. }\frac{1}{2} \log \langle w, w \rangle_{H_T} < - \epsilon_1 \log X \}}{ \#\mathcal{F}(X)} = 0. \]
This follows from the equidistribution of $\mathcal{R}$ on irreducible integral orbits, more specifically Corollary \ref{corollary_polyswithsmallnormvectorrare}.
\end{proof}
\bibliographystyle{alpha}
\bibliography{references}

\end{document}